\documentclass[12pt]{article}
 
\newif\ifSHOWEXTRA
\SHOWEXTRAtrue
\SHOWEXTRAfalse

\usepackage[
margin=1in]{geometry}
\usepackage{amsmath,amsthm,amssymb, latexsym,amsfonts,bbm, scrextend}
\usepackage[usenames,dvipsnames]{xcolor}
\usepackage{fancyhdr}
\usepackage{array}
\usepackage {tikz}
\usetikzlibrary {positioning}
\DeclareFontFamily{U}{rcjhbltx}{}
\DeclareFontShape{U}{rcjhbltx}{m}{n}{<->rcjhbltx}{}
\DeclareSymbolFont{hebrewletters}{U}{rcjhbltx}{m}{n}

\let\aleph\relax\let\beth\relax
\let\gimel\relax\let\daleth\relax

\DeclareMathSymbol{\aleph}{\mathord}{hebrewletters}{39}
\DeclareMathSymbol{\beth}{\mathord}{hebrewletters}{98}
\DeclareMathSymbol{\gimel}{\mathord}{hebrewletters}{103}
\DeclareMathSymbol{\daleth}{\mathord}{hebrewletters}{100}

\DeclareMathSymbol{\lamed}{\mathord}{hebrewletters}{108}
\DeclareMathSymbol{\mem}{\mathord}{hebrewletters}{109}
\DeclareMathSymbol{\ayin}{\mathord}{hebrewletters}{96}
\DeclareMathSymbol{\tsadi}{\mathord}{hebrewletters}{118}
\DeclareMathSymbol{\qof}{\mathord}{hebrewletters}{113}
\DeclareMathSymbol{\samech}{\mathord}{hebrewletters}{115}
\DeclareMathSymbol{\shin}{\mathord}{hebrewletters}{152}

\definecolor {processblue}{cmyk}{0.96,0,0,0}
\definecolor {processred}{cmyk}{0,1,.5,.25}

\newcommand{\given}{\hspace{3pt}\vline\hspace{3pt}}
\newcommand{\givenn}{\hspace{3pt}|\hspace{3pt}}
\newcommand{\Pm}{{\bf P}}

\newcommand{\pre}{\upsilon}
\newcommand{\PRE}{\Upsilon}
\newcommand{\bOmega}{\boldsymbol{\Omega}}
\newcommand{\bomega}{\boldsymbol{\omega}}
\newcommand{\tomega}{\tilde{\bomega}}
\newcommand{\BigSpace}{\text{\scalebox{1.5}{$\samech$}}}
\newcommand{\tuple}{\shin}
\newcommand{\siteb}{b}
\newcommand{\setM}{\mathcal{S}}
\newcommand{\Em}{{\bf E}}
\newcommand{\R}{\mathbb R}
\newcommand{\Pa}{\mathbb P}

\newcommand{\Ea}{\mathbb E}

\newcommand{\N}{\mathbb N}

\newcommand{\Z}{\mathbb Z}

\newcommand{\vel}{\overset{\star}{v}}
\newcommand{\vect}[1]{\hat{#1}}
\usepackage{tikz}
\usepackage{pgfplots}
\usepackage{amsmath}
\usepackage{eqnarray}
\newcounter{CommentCounter}

\usepackage[mathscr]{euscript}
 \let\mathscr\relax
\usepackage[scr]{rsfso}
\usepackage{amsthm}
\usepackage{amssymb}
\usepackage{multicol}
\usepackage{enumerate}
\usepackage[colorlinks=true, pdfstartview=FitV, linkcolor=blue,
citecolor=blue, urlcolor=blue]{hyperref}

\newcommand{\nodetype}[6]{\begin {tikzpicture}[-latex ,auto ,node distance =2 cm and 2cm ,on grid ,
semithick ,
state/.style ={ circle ,top color =white , bottom color = #5!20 ,
draw,#5 , text=#5 , minimum width =1 cm}]
\node[state] (C)            
{#6};
\node[draw=none] (A) [above =of C] {};
\node[draw=none] (B) [right =of C] {};
\node[draw=none] (D) [left =of C] {};
\node[draw=none] (E) [below =of C] {};
\path (C) edge node[left=.01 cm]{#1} (A);
\path (C) edge node[below=.01 cm]{#2} (D);
\path (C) edge node[right=.01 cm]{#3} (E);
\path (C) edge node[above=.01 cm]{#4} (B);
\end{tikzpicture}}

\newcommand{\one}[1]{\mathbbm{1}_{\{\text{#1}\}}}
\newtheorem{claim}{Claim}
 
\newtheorem*{rem}{Remark}
\newtheorem{thm}{Theorem}
\newtheorem{lem}{Lemma}

\newtheorem{prop}{Proposition}
\newtheorem*{prop*}{Proposition}

\newtheorem{cor}[thm]{Corollary}

\newtheorem{exmp}{Example}
\theoremstyle{definition}

\newtheorem{question}{Question}
\pgfplotsset{compat=1.12}

\begin{document}

\title{Random Walks in Random Environments with Rare Anomalies}
\author{Daniel J. Slonim}

\date{\today}
\maketitle
\abstract{We study random walks in i.i.d. random environments on $\Z^d$ when there are two basic types of vertices, which we call ``blue'' and ``red''. Each color represents a different probability distribution on transition probability vectors.
We introduce a method of studying these walks that compares the expected amount of time spent at a specific site on the event that the site is red with the expected amount of time spent there on the event that the site is blue. This method produces explicit bounds on the asymptotic velocity of the walk. We recover an early result of Kalikow, but with new  bounds on the velocity. Next, we consider a ``rare anomaly'' model where the vast majority of sites are blue, and blue sites are uniformly elliptic, with some almost-sure bounds on the quenched drift. We show that if the red sites satisfy a certain uniform ellipticity assumption in two fixed, non-parallel directions, then even if red sites break the almost-sure bounds on the quenched drift, making red sites unlikely enough lets us obtain bounds on the asymptotic velocity of the walk arbitrarily close to the bounds on the quenced drift at blue sites. Significantly, the required proportion $p^*$ of blue sites to do this does not depend on the distribution of red sites, except through the uniform ellipticity assumption in two directions. Our proof is based on a coupling technique, where two walks run in environments that are the same everywhere except at one vertex. They decouple when they hit that vertex, and our proof is driven by bounds on how long it takes to recouple. We then demonstrate the importance of the i.i.d. assumption by providing a counterexample to the statement of the theorem with this assumption removed. We conclude with open questions.
    
    \medskip\noindent {\it MSC 2020:}
    60G50 
    60J10 
    60K37 
    \\
    {\it Kewords:}
    random walk,
    random environment,
    transience,
    ballistic,
    asymptotic velocity,
    continuity,
    rare anomalies
}


\section{Introduction}

This paper provides continuity-type results for asymptotic velocity of random walks in random environments (RWRE) on $\Z^d$ for $d\geq2$, giving  explicit velocity bounds in terms of the environment at only one site. The model involves a random Markov chain with $\Z^d$ as the state space. A typical assumption is that nearest-neighbor transition probability vectors at individual states are drawn in an i.i.d. way according to some underlying distribution. For instance, in \cite{Kalikow1981}, S. Kalikow presented the following i.i.d. \textit{two-vertex} model: Fix $p\in(0,1)$ and positive probability vectors $<b_1,b_2,b_3,b_4>$ and $<r_1,r_2,r_3,r_4>$. Say that all sites in $\Z^2$ are independently ``blue'' with probability $p$ and ``red'' with probability $1-p$. 
Then, letting $e_1,e_2,e_3,e_4$ respectively denote the ``right'', ``up'', ``left'', and ``down'' axis unit vectors, say that a walk ${\bf X}=(X_n)_{n=0}^{\infty}$ steps in direction $e_i$ from a blue site with probability $b_i$, and from a red site with probability $r_i$.

\begin{figure}[ht]
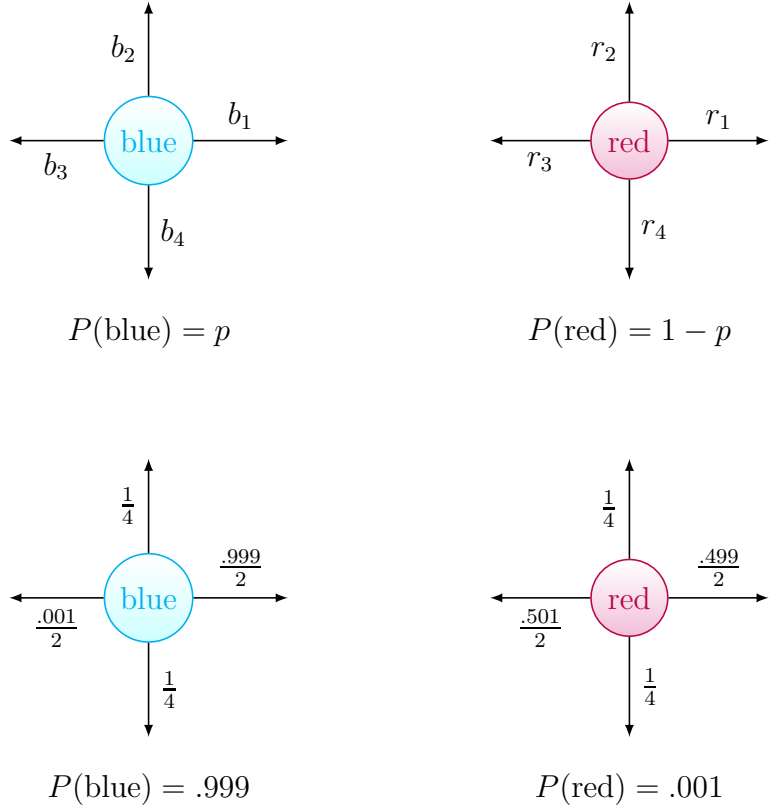

    \centering
    \begin{tabular}{ccc}
     \nodetype{$b_2$}{$b_3$}{$b_4$}{$b_1$}{processblue}{blue}    
     & \hspace{.5in}
     &
    \nodetype{$r_2$}{$r_3$}{$r_4$}{$r_1$}{processred}{red}
     \\\notag
     $P(\text{blue})=p$  &  & $P(\text{red})=1-p$
      \\\notag
  \vspace{.25in}
    \\\notag
    \nodetype{$\frac14$}{$\frac{.001}{2}$}{$\frac14$}{$\frac{.999}{2}$}{processblue}{blue}&
    \hspace{.5in}&
    \nodetype{$\frac14$}{$\frac{.501}{2}$}{$\frac14$}{$\frac{.499}{2}$}{processred}{red}
    \\\notag
    $P(\text{blue})=.999$  &  & $P(\text{red})=.001$
    \end{tabular}
    \caption{Kalikow's two-vertex model is depicted on the top. On the bottom is depicted the two-vertex model with the specific values Kalikow presented at the beginnng of \cite{Kalikow1981}.}
    \label{fig:2Vertex}
\end{figure} 

When $d\geq2$, RWRE on $\Z^d$ are far harder to understand than when $d=1$. To illustrate this difficulty, the paper \cite{Kalikow1981} opens rather provocatively with the following question about the two-vertex model on $\Z^2$: If $r_2=r_4=b_2=b_4=\frac14$ (so that there is deterministically no vertical drift at any site and all vertical step probabilities are moderate), if $b_1=\frac{.999}{2}$, $b_3=\frac{.001}{2}$, $r_1=\frac{.499}{2}$, if $r_3=\frac{.501}{2}$ (so that there is a very strong rightward drift at blue sites and a very mild leftward drift at red sites), and if $p=.999$ (so that sites are overwhelmingly likely to be blue), then can one prove that the walk is almost surely transient to the right?
Although an affirmative answer seems obvious at first, \cite{Kalikow1981} takes a good deal of work to provide a proof, which comes as a result of showing that the model satisfies a more general but non-explicit criterion known as ``Kalikow's criterion''. A sufficient condition for Kalikow's criterion comprises the hypothesis of the following theorem.


\begin{thm}\label{thm:RecoveringKalikow+Bounds}
    Suppose $b_1>b_3$ and $r_3>r_1$. Moreover, suppose
    \begin{equation}\label{eqn:220}
        \frac{p(b_1-b_3)}{(1-p)(r_3-r_1)}>M:=\max\left\{\frac{b_i}{r_i};i=1,2,3,4\right\}.
    \end{equation}
    Then 
    \begin{enumerate}[(a)]
        \item the walk is transient to the right---that is, $\lim_{n\to\infty}X_n\cdot e_1=\infty$, almost surely;
        \item
     The walk is ballistic; in fact, with probability 1,
    \begin{equation}\label{eqn:Bounds}
    \hspace{-.25in}
        \frac{p(b_1-b_3)-(1-p)M(r_3-r_1)}{p+(1-p)M}
        \leq
        \lim_{n\to\infty}\frac{X_n\cdot e_1}{n}
        \leq
         \frac{b_1-b_3}{1+(1-p)[(b_1-b_3)+(r_3-r_1)]}.
         \hspace{.25in}
    \end{equation}
    \end{enumerate}
\end{thm}
Part (a) is proven in \cite{Kalikow1981} as the Corollary to Theorem 1, but we give a new proof. Sznitman and Zerner in \cite{Sznitman&Zerner1999} showed that under Kalikow's criterion, the walk is ballistic (has nonzero limiting speed). The paper \cite{Zerner2000} expressed the asymptotic velocity of models satisfying Kalikow's condition in terms of Lyapunov exponents associated with the walk. The explicit lower bound in \eqref{eqn:Bounds} is new, and comes as a result of our new proof of Part (a). For an upper bound, the maximum rightward quenched drift $b_1-b_3$ works, but the fact that a ballistic walk must visit ``fresh'' sites---where the annealed drift controls---at regular intervals yields a slight improvement.


To understand some implications and limits of Theorem \ref{thm:RecoveringKalikow+Bounds}, consider the following model, where all vertical steps deterministically have probability $\frac14$ and blue sites have a moderate drift of $.1$ to the right, while at a red site, the probability of stepping to the right is $\delta$.
\begin{figure}[ht]
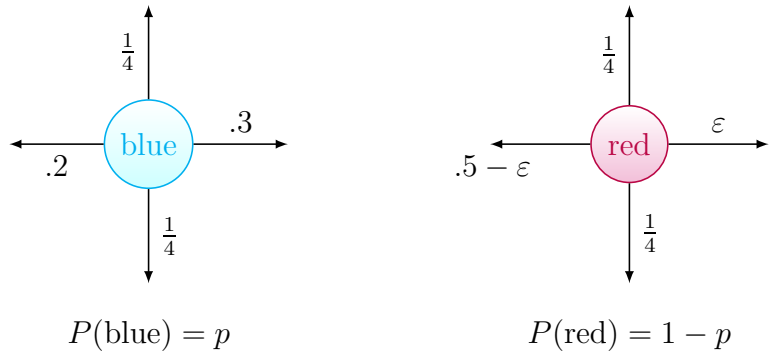

    \centering
    \begin{tabular}{ccc}
     \nodetype{$\frac14$}{$.2$}{$\frac14$}{$.3$}{processblue}{blue}    
     & \hspace{.5in}
     &
    \nodetype{$\frac14$}{\hspace{-.5in}$.5-\delta$}{$\frac14$}{$\delta$}{processred}{red}
     \\\notag
     $P(\text{blue})=p$  &  & $P(\text{red})=1-p$
    \end{tabular}
    \caption{The case described in Theorem \ref{thm:SpecialCase}}
    \label{fig:SpecialCase}
\end{figure} 
Here, if $\delta$ is small, then the $M$ in \eqref{eqn:220} is large, requiring $p$ to be close to 1 in order for the model to satisfy the hypothesis of Theorem \ref{thm:RecoveringKalikow+Bounds}. 
Indeed, asymptotically, $p$ must be at least $1-\frac23\delta$.  And if $\delta=0$, then Theorem 1 can never show transience to the right no matter how close $p$ comes to 1.

A goal of the present paper is to show transience to the right and ballisticity under milder assumptions that do not require $p$ to be taken arbitrarily close to 1. Indeed, we show that for sufficiently large (but fixed) $p$, $\delta$ can be taken as small as one wishes---including 0---and transience to the right will hold, with a positive lower bound on the velocity that is uniform in $\delta$ and that can be made arbitrarily close to 0.1, the strength of the rightward drift at blue sites.

\begin{thm}\label{thm:SpecialCase}
    Let $\delta\geq0$. Consider the i.i.d. two-vertex model where $b_2=b_4=r_2=r_4=\frac14$, $b_1=.3$, $b_3=.2$, $r_1=\delta$, and $r_3=\frac12-\delta$. For every $\varepsilon>0$, there exists a $p^*\in(0,1)$ such that whenever $p>p^*$, the almost-sure limiting speed $\lim_{n\to\infty}\frac{X_n\cdot e_1}{n}$ is at least $0.1-\varepsilon$.
\end{thm}

\begin{rem}
    An analogous result would be false in one dimension. If blue sites send the walker to the right with probability $b_1>\frac12$ and to the left with probability $b_2=1-b_1$, and red sites send the walker to the right with probability $\delta$ and left with probability $1-\delta$, then the results of \cite{Solomon1975} imply that for any fixed $b_1,p\in(0,1)$, the walk can be made transient to the left by taking $\delta$ sufficiently small. However, this is due to the fact that a walk that is transient to the right must go through every vertex; there is no way to ``walk around'' a bad vertex. In higher dimensions, this is no longer the case.
\end{rem}

We prove Theorem \ref{thm:SpecialCase} as a special case of the more general Theorem \ref{thm:MainThm}, which requires more notation to state. We call the model for this theorem the \textit{rare anomaly model}. The blue sites are allowed to be drawn from one non-degenerate distribution with almost-sure bounds on the quenched drift at a given site, and the red sites---the rare anomalies---are allowed to be drawn from another, potentially very different distribution that breaks these bounds. The theorem will say that if red sites are sufficiently rare, the limiting velocity of the walk can be given bounds arbitrarily close to the bounds on the quenched drift at blue sites. 


The rare anomaly model is natural to think about in real-world applications of simple random walks (i.e., random walks in homogeneous environments). Indeed, for any model where a (perhaps biased) simple random walk is used to model physical movement in a homogeneous medium, the assumption of perfect homogeneity is likely to be unrealistic in the real world. Nevertheless, in many cases it may be reasonable to assume that the behavior of the walk in an environment with anomalies---if they are sufficiently rare---will be similar to the behavior of a walk in a perfectly homogeneous environment, and that motion in a nearly homogeneous medium can therefore be reasonably modeled by a random walk in a perfectly homogeneous environment. Our work here makes this idea explicit for the first-order behavior of certain RWRE models. 

In the case where the support of $\mu_b$ is either a single vector or contained in a neighborhood of a single vector, 
the rare anomaly model is a perturbation---at the level of distributions---of a simple random walk. Perturbations of simple random walks have been studied in papers such as \cite{Bolthausen&Sznitman2002}, \cite{Sznitman2003}, \cite{Sabot2004}, and \cite{Bolthausen&Zeitouni2007}. However, these focus on distributions with support contained entirely in small neighborhoods of a fixed vector, so that the \textit{transition probability vector} at every site is almost surely a small perturbation of a fixed vector. Our model allows us to study \textit{distributions} that are perturbations of a degenerate distribution supported on a single vector. To make this idea explicit, we state and prove Corollary \ref{cor:continuity}, which considers a large subset of the $2d$-simplex, and says that the almost-sure limiting velocity of the walk\footnote{Actually, the object we consider is not the limiting velocity itself but a probability measure consisting of either a Dirac mass at the limiting velocity (if it is deterministic), or half a Dirac mass at each of two possible limiting velocities. This technicality is due to the fact that we cannot yet rule out the possibility that there are two possible limiting velocities.}, as a function of the probability distribution on that subset of the simplex, is continuous with respect to the topology of weak convergence of distributions. 

Since Kalikow's paper \cite{Kalikow1981}, a good deal of work has been done studying ballisticity of RWRE in dimensions higher than 1. In the special case of Dirichlet environments, ballisticity has been completely characterized for dimensions 3 and higher in terms of parameters that are very simple to read off from the distribution of the environment at a single site \cite{Sabot2013}. For general distributions, several abstract sufficient conditions for ballisticity have been given, and the relationships between these conditions has been a highly active area of study; see \cite{Sznitman2001}, \cite{Drewitz&Ramirez2011}, \cite{Berger&Drewitz&Ramirez2014}, \cite{Guerra&Ramirez2019}. Some of these abstract conditions can, in theory, be directly checked; such conditions are dubbed ``effective criteria'' (see, e.g., \cite{Sznitman2002} \cite{Berger&Drewitz&Ramirez2014}). Such criteria often involve such heavy computation that even a computer cannot check them for most specific examples; however, they can be used to prove existence of distributions on environments satisfying certain ballisticity properties.

To the author's knowledge, Theorems \ref{thm:RecoveringKalikow+Bounds} and \ref{thm:MainThm} are the first results that give explicit, non-trivial bounds on the limiting velocity for a class of i.i.d., ballistic RWRE in dimension higher than 1 that can be directly computed in terms of the distribution of the environment at a single site.

We formally introduce our model and state Theorem \ref{thm:MainThm} in subsection \ref{subsec:Model&Statement}. In subsection \ref{subsec:Continuity}, we discuss Corollary \ref{cor:continuity}, which justifies the idea that Theorem \ref{thm:MainThm} is a continuity-type result. We describe the main proof ideas in Subsection \ref{subsec:ProofIdeas}. Section \ref{sec:Prelim} describes a weakened version of an ellipticity assumption and reviews a previously known law of large numbers under this assumption that will be crucial for our results.
Section \ref{sec:Criterion} gives a criterion for directional transience and ballisticity, which will be used both for Theorem \ref{thm:MainThm} and Theorem \ref{thm:RecoveringKalikow+Bounds}. Section \ref{sec:Kalikow} proves Theorem \ref{thm:RecoveringKalikow+Bounds}. Section \ref{sec:RareObs} proves  Theorem \ref{thm:MainThm}. Section \ref{sec:Counterexample} demonstrates the importance of the i.i.d. assumption in Theorem \ref{thm:MainThm} by providing a non-i.i.d. counterexample with good mixing properties  where, for any $p\in(0,1)$, blue and red sites are as described in Theorem \ref{thm:SpecialCase} and sites are blue with probability $p$, but where the walk does not have positive speed in direction $e_1$. Section \ref{sec:Qs} concludes by discussing some open questions.

\subsection{Rare anomaly model and main result statement}\label{subsec:Model&Statement}

Let $\mathcal{N}=\{e_1,\ldots,e_{2d}\}$ be the set of nearest neighbors of the origin in $\Z^d$, with $e_{i+d}=-e_i$, $1\leq i\leq d$. Define 
\[\Sigma_{\mathcal{N}}:=\left\{\xi:\mathcal{N}\to[0,1]:\sum_{i=1}^{2d}\xi(e_i)=1\right\}.\]
This is the simplex---or the set of probability vectors---on $\mathcal{N}$. Let 
\begin{equation}\notag\Omega:=\prod_{x\in \Z^d}\Sigma_{\mathcal{N}}:=\{(\omega^x)_{x\in\Z^d}:\omega^x\in\Sigma_{\mathcal{N}}\text{ for all }x\in\Z^d\}.\end{equation} We refer to an element $\omega=(\omega^x)_{x\in\Z^d}$ as an \textit{environment} on $\Z^d$. 
If $X\subset \Z^d$, then let $\omega^X:=(\omega^x)_{x\in X}$.

For a given environment $\omega$ and $x\in \Z^d$, we can define $P_{\omega}^x$ to be the measure on $(\Z^d)^{\N}$ (with the natural sigma field) giving the law of a Markov chain $(X_n)_{n=0}^{\infty}$ started at $x$ with transition probabilities given by $\omega$. That is, $P_{\omega}^x(X_0=x)=1$, and for $n\geq1$, $y\in V$, $P_{\omega}^x(X_{n+1}-X_n=y\given X_0,\ldots,X_n)=\omega^{X_n}(y)$ almost surely.

Let us generalize Kalikow's two-vertex model by allowing for two \textit{classes} of vertices. Let ``blue'' vertices have transition probability vectors distributed according to one measure $\mu_b$, and ``red'' vertices have transition probability vectors distributed according to another measure $\mu_r$. In order to keep track of which sites are blue and red, we will describe the notion of a \textit{tagged environment}, which keeps track of which sites are blue and which are red as well as their transition probability vectors.

Let
\begin{equation}\notag\PRE:=\{\text{blue},\text{red}\}^{\Z^d}:=\{(\pre(x))_{x\in\Z^d}:\pre(x)\in\{\text{blue},\text{red}\}\text{ for all }x\in\Z^d\}.\end{equation}
be the set of ``pre-environments'', and let \begin{equation}\notag\bOmega:=\PRE\times\Omega\end{equation}
be the set of ``tagged environments''---that is, environments where each site is assigned a color as well as a set of transition probabilities.\footnote{Because we are taking an i.i.d. assumption, we could assume without loss of generality that the measures $\mu_b$ and $\mu_r$ are mutually singular. This would allow us to make the tagged environment a function of the environment, which would avoid the need to define $\PRE$ and $\bOmega$. However, defining the tagged environment at this stage lets us ask natural questions about what happens if we slightly relax the i.i.d. assumption; see for instance Question \ref{quest:RelaxIID}.} 
We use $\bomega$ to denote a generic element $(\pre,\omega)\in\bOmega$.

Now if we give $\Sigma_{\mathcal{N}}$ the subspace  topology, $\{\text{blue},\text{red}\}$ the discrete topology, and $\Omega$, $\PRE$, and $\bOmega$ the product topologies, we can let $\boldsymbol{\mathcal{F}}$ be the Borel sigma field on $\bOmega$.  
Let $\Pm$ be a probability measure on $(\bOmega,\boldsymbol{\mathcal{F}})$. For a given $x\in \Z^d$, we let $\Pa^x=\Pm\times P_{\omega}^x$ be the measure on $\bOmega\times (\Z^d)^{\N}$ induced by both $\Pm$ and $P_{\omega}$. That is, for measurable events $A\subset\bOmega,B\subset (\Z^d)^{\N}$,
\begin{equation}\label{eqn:362}
\Pa^x(A\times B)=\int_{A}P_{\omega}^x(B)\Pm(d\bomega)
\end{equation}
In particular, $\Pa^x(\bOmega\times B)=\Em\left[P_{\omega}^x(B)\right]$. 
For convenience, we commit a small abuse of notation by using $\Pa^x$ to refer both to the measure we've described on $\bOmega \times (\Z^d)^{\N}$ and also to its marginal $\Pa^x(\bOmega \times \cdot)$ on $(\Z^d)^{\N}$. 
We call a measure $P_{\omega}^x$ on $(\Z^d)^{\N}$ a {\em quenched measure} of a random walk in random environment on $\Z$ started at $x$, and we call the measure $\Pa^x$ 
the {\em annealed measure}.


 Let $\mu_b$ and $\mu_r$ be two probability measures on $\Sigma_{\mathcal{N}}$ with the following properties.

\begin{enumerate}
    \item $\mu_b$ satisfies a \textit{uniform ellipticity} assumption; there exists $\kappa>0$ such that every vector in the support of $\mu_b$ has all $2d$ components at least $\kappa$.
    \item $\mu_r$ satisfies the following ``two-direction'' uniform ellipticity assumption: there exists $\kappa>0$ and two distinct values $\eta,\theta\in\{1,\ldots, d\}$ with $|\eta-\theta|\notin\{0,d\}$ (i.e., the vectors associated with $e_{\eta}$ and $e_{\theta}$ are orthogonal) such that every vector in the support of $\mu_r$ has $e_{\eta}$- and $e_{\theta}$-components at least $\kappa$.
\end{enumerate}

Now let $p\in(0,1)$. Let $\Pm$ be the measure on $\bOmega$ such that 
\begin{enumerate}
    \item The pairs $(\pre(x),\omega^x)$ are i.i.d.;
    \item For each $x\in\Z^d$, $\pre(x)=\text{blue}$ with probability $p$, and $\pre(x)=\text{red}$ with probability $1-p$.
    \item Conditioned on $\pre(x)=\text{blue}$, the law of $\omega^x$ is $\mu_b$;
    \item Conditioned on $\pre(x)=\text{red}$, the law of $\omega^x$ is $\mu_r$.
\end{enumerate}

Note that Kalikow's two-vertex model is simply the special case where $d=2$ and where $\mu_b$ and $\mu_r$ are each a Dirac mass at one probability vector with all positive entries.

We are now able to state our main result, a more general version of Theorem \ref{thm:SpecialCase}.

\begin{thm}\label{thm:MainThm}
 Let $S^{d-1}$ denote the unit sphere in $\R^d$, and let $\vec{u}\in S^{d-1}$. Suppose $v_1^b,v_2^b$ are such that every vector $\xi$ in the support of  $\mu_b$ has $\sum_{i=1}^{2d}\xi(e_i)e_i\cdot\vec{u}\in[v_1^b,v_2^b]$. For every $\varepsilon>0$, there exists $p^*<1$, depending only on $\varepsilon$, $d$, $v_1^b$, and $\kappa$, such that if $p>p^*$, then $\Pa^0$--almost surely, 
 \begin{equation}
  v_1^b-\varepsilon\leq \lim_{n\to\infty}\frac{X_n\cdot \vec{u}}{n}\leq v_2^b+\varepsilon.
 \end{equation}
\end{thm}

\begin{rem}\label{rem:exactbound}
    One can see from the proofs that if $\varepsilon<1$, then it will suffice\footnote{From the proofs, one will also see that $p^*$ must be greater than $p_c(d)$, the critical site percolation threshold for $\Z^d$. However, this threshold can only decrease with dimension, and for $d=2$, one can easily see from, e.g., \cite[Theorem 1.33 and Theorem 11.11]{Grimmett1999}, that the critical threshold is no more than $\frac{15}{16}$, which is less than $1-\frac{\kappa^6\varepsilon}{13d}$ for any non-trivial $\varepsilon$.}  to take $p^*=1-\frac{\kappa^6\varepsilon}{13d}$. Thus, the result can be used to obtain explicit bounds on the limiting velocity. 
\end{rem}

Note that in particular the value of $p^*$ depends on the distribution $\mu_r$ only through $\kappa$. In other words, one can get a lower bound on the speed of the walk in direction $e_1$ that is arbitrarily close to $v_1^b$ simply by controlling how rare the exceptional sites are, without regard for what the exceptional sites look like, aside from the uniform ellipticity in two of the $2d$ directions.

\begin{exmp}
    In Theorem \ref{thm:SpecialCase}, let $p=1-10^{-8}$. Then, using Remark \ref{rem:exactbound} (with $\kappa=0.2$), we can see that $\lim_{n\to\infty}\frac{X_n\cdot e_1}{n}\geq0.095$, regardless of the value of $\delta$. 
\end{exmp}

\begin{exmp}[Nearly non-nestling]
    Let $\mu_b$ and $\mu_r$ be as in the statement of Theorem \ref{thm:MainThm}, and let $\mu_b$ satisfy the following \emph{non-nestling} assumption: there is a unit vector $\vec{u}\in S^{d-1}$ and a $v_1>0$ such that for every $\xi$ in the support of $\mu_b$, $\sum_{i=1}^{2d}\xi(e_i)e_i\cdot\vec{u}\geq v_1.$ Then a RWRE where transition probability vectors are drawn according to $\mu_b$ at every site is of course transient in direction $\vec{u}$ with asymptotic speed $\lim_{n\to\infty}\frac{X_n\cdot\vec{u}}{n}$ at least $v_1$. By Theorem \ref{thm:MainThm}, even if the existence of red sites causes the overall distribution to break the non-nestling assumption, the walk is still almost-surely transient in direction $\vec{u}$ with a lower bound on the asymptotic velocity that can be made arbitrarily close to $v_1$, provided $p$ is sufficiently close to 1. This example can be seen as a generalization of the non-nestling regime from \cite{Zerner1998}, \cite{Komorowski&Krupa2003}, \cite{Rassoul-Agha&Seppalainen2009}.

\begin{exmp}[Double perturbation of a zero-drift simple random walk]
    Let $d=4$, and let $\varepsilon>0$. Let $\mu_b$ be a probability measure on $\Sigma_{\mathcal{N}}$ such that every vector in the support of $\mu_b$ has each component within $\frac{\varepsilon}{4}$ of $\frac18$. Thus, under the measure on environments associated to $\mu_b$, every vector is a small perturbation of a vector assigning equal probability to all 8 directions, and the magnitude of any quenched drift is at most $\varepsilon$. In general, it is not known whether under the annealed measure corresponding to $\mu_b$, the limiting velocity is deterministic, but nevertheless using almost-sure bounds on the quenched drift, the limiting velocity must almost surely have magnitude no more than $\varepsilon$. Now let $\mu_r$ be any probability measure on $\Sigma_{\mathcal{N}}$ under which the probabilities associated to directions $e_1$ and $e_2$ are almost surely at least $0.1$. If transition probability vectors are drawn from $\mu_b$ with probability $p<1$ and from $\mu_r$ with probability $1-p$, then we no longer have almost-sure bounds on quenched drift (except in two directions, and not very good ones at that). Nevertheless, by Theorem \ref{thm:MainThm} and Remark \ref{rem:exactbound}, we can show that if $1-p<\frac{\varepsilon}{208,000,000}$, then the limiting velocity in any direction $e_1$, $e_2$, $e_3$, $e_4$ has magnitude at most $\varepsilon$, so the overall limiting velocity almost surely has magnitude at most $2\varepsilon$.
\end{exmp}

\end{exmp}

\subsection{A continuity result}\label{subsec:Continuity}

Theorem \ref{thm:MainThm} is a continuity-type result in the sense that it shows that perturbing the distribution on environments slightly does not change the asymptotic velocity ``too much''. Here, we give a corollary to Theorem \ref{thm:MainThm} that is stated directly in terms of continuity, partially extending a known continuity result from the one-dimensional case.

The paper \cite{Solomon1975} showed that in the i.i.d. case for $d=1$, when the law on environments is $\bigotimes_{z\in\Z}\mu$ for some probability measure $\mu$ on $\Sigma_{\mathcal{N}}$ the walk has a deterministic asymptotic limiting velocity
$v:=\lim_{n\to\infty}\frac{X_n}{n}$ given by

\begin{equation}\label{eqn:236}
    v=\begin{cases}
        \frac{1-\Em\left[\frac{\omega^0(-e_1)}{\omega^0(e_1)}\right]}{1+\Em\left[\frac{\omega^0(-e_1)}{\omega^0(e_1)}\right]} &\text{if }-\Em\left[\log\left(\frac{\omega^0(-e_1)}{\omega^0(e_1)}\right)\right]>0, \Em\left[\frac{\omega^0(-e_1)}{\omega^0(e_1)}\right]<1;
        \\
        -\frac{1-\Em\left[\frac{\omega^0(e_1)}{\omega^0(-e_1)}\right]}{1+\Em\left[\frac{\omega^0(e_1)}{\omega^0(-e_1)}\right]} &\text{if }-\Em\left[\log\left(\frac{\omega^0(-e_1)}{\omega^0(e_1)}\right)\right]<0, \Em\left[\frac{\omega^0(e_1)}{\omega^0(-e_1)}\right]<1;
        \\
        0&\text{otherwise,}
    \end{cases}
\end{equation}
provided the expectation $\Em\left[\log\left(\frac{\omega^0(-e_1)}{\omega^0(e_1)}\right)\right]$ is well defined. From results of \cite{Key1984}, \cite{Sznitman&Zerner1999}, and \cite{Zerner2002}, one can see that a deterministic limiting velocity $v$ must exist even when the expectation is not well defined.

The limiting velocity $v$ may be viewed as a function of the probability measure $\mu$ on  $\Sigma_{\mathcal{N}}$, and it is natural to ask about the continuity of such a function. Let $\mathcal{M}_1\left(\Sigma_{\mathcal{N}}\right)$ be the set of probaiblity measures on $\Sigma_{\mathcal{N}}$, endowed with the topology of weak convergence. Is $v$ continuous as a function on this space? The initial answer is no: if $\mu$ is a measure with $v(\mu)>0$, one may take a sequence $\mu_n$ converging to $\mu$, where each $\mu_n$ is a convex combination of $\mu$ and a Dirac mass at the probability vector $\xi$ with $\xi(e_1)=0$ and $\xi(-e_1)=1$. Under the annealed measure associated to each $\mu_n$, the walk cannot be so much as transient to the right, since sites will exist where the walk gets sent to the left on every visit, and such sites cannot be passed. Therefore, $v(\mu_n)\leq0$ for all $n$ (one can show that eventually $v(\mu_n)=0$).

Nevertheless, if we restrict our attention to a subset of $\Sigma_{\mathcal{N}}$ where a \textit{uniform ellipticity} assumption is satisfied, then the answer to the continuity question becomes yes. For $\kappa>0$, 
define 
$$
\Sigma_{\mathcal{N}}^{\kappa}:=\{\xi\in\Sigma_{\mathcal{N}}:\xi(e_1),~\xi(e_2)\geq\kappa\}.
$$
If the support of $\mu$ is contained in $\Sigma_{\mathcal{N}}^{\kappa}$, then the product measure $\Pm$ is \textit{uniformly elliptic} with ellipticity constant $\kappa$. Now if $\mathcal{M}_1\left(\Sigma_{\mathcal{N}}^{\kappa}\right)$ is the set of probability measures on $\Sigma_{\mathcal{N}}^{\kappa}$, endowed with the topology of weak convergence, then the limiting velocity $v$ is always given by \eqref{eqn:236}, since the well-definedness of $\Em\left[\log\left(\frac{\omega^0(-e_1)}{\omega^0(e_1)}\right)\right]$ is automatic, and one can easily check that the function $v$ is continuous on $\mathcal{M}_1\left(\Sigma_{\mathcal{N}}^{\kappa}\right)$ using bounded convergence. We thus have

\begin{prop*}
Let $d=1$. The function $v$ is continuous on $\Sigma_{\mathcal{N}}^{\kappa}$ for every $\kappa>0$.     
\end{prop*}

We now consider $d\geq2$. Here there is no formula for a limiting velocity $v$, and for $d\geq3$, it is not even known whether a deterministic limiting velocity exists. However, by the results of \cite{Sznitman&Zerner1999} and \cite{Zerner2002}, it is known that a (possibly random) limiting velocity exists and can be supported on at most two vectors, which must be negative multiples of each other (in dimension 2 the limiting velocity is known to be deterministic by the results of \cite{Zerner&Merkl2001}). 

Rather than restricting all the way to the uniformly elliptic setting $\Sigma_{\mathcal{N}}^{\kappa}$; we only need to assume that uniform ellipticity is satisfied in two directions. Fix directions $\eta,\theta\in\{1,\ldots, d\}$ with $|\eta-\theta|\notin\{0,d\}$ and $\kappa>0$.  Define
\[\Sigma_{\mathcal{N}}^{\kappa,\eta,\theta}=\left\{\xi\in\Sigma_{\mathcal{N}}:\sum_{i=1}^{2d}\xi(e_i)=1,~e_{\theta},e_{\eta}\geq\kappa\right\}.\]
This is the part of the simplex $\Sigma_{\mathcal{N}}$ where $e_{\theta},e_{\eta}\geq\kappa$. It is naturally embedded as a compact subset of $\R^{2d-1}$. Endow it with the subspace topology, and let $\mathcal{M}_1(\Sigma_{\mathcal{N}}^{\kappa,\eta,\theta})$ be the set of Borel probability measures on this set, with the topology of weak convergence. For a measure $\mu\in\mathcal{M}_1(\Sigma_{\mathcal{N}}^{\kappa,\eta,\theta})$, consider a RWRE where transition probabilities are drawn in an i.i.d. way at all sites according to $\mu$. Let $B_1^1(0)$ be the $L^1$ unit ball centered at the origin in $\R^d$. A nearest-neighbor walk must have $\frac{X_n}{n}\in B_1^1(0)$ for all $n$, so the asymptotic velocity must be in $B_1^1(0)$. Let $\mathcal{M}_1(B_1^1(0))$ be the space of probability measures on $B_1^1(0)$, with the topology of weak convergence. Define the function $\vel:\mathcal{M}_1(\Sigma_{\mathcal{N}}^{\kappa,\eta,\theta})\to\mathcal{M}_1(B_1^1(0))$ by 
\begin{equation}
    \vel(\mu)=
    \begin{cases}
    \delta_{\vec{v}} &\parbox{4in}{if, when sites are drawn in an i.i.d. way according to $\mu$, the limiting velocity is almost surely $\vec{v}\in\R^d$};
    \\
    \vspace{1pt}
    \\
    \frac12\delta_{\vec{v}}+\frac12\delta_{\vec{w}}&\parbox{4in}{if, when sites are drawn in an i.i.d. way according to $\mu$, the limiting velocity is almost surely either $\vec{v}\in\R^d$ or $\vec{w}\in\R^d$, each with positive probability.}
    \end{cases}
\end{equation}
Thus, $\vel(\mu)$ is a Dirac mass at the almost-sure limiting velocity of a RWRE with sites drawn in an i.i.d. way according to $\mu$, if the limiting velocity is deterministic, and half a Dirac mass at each point if the limiting velocity is supported on two points.\footnote{Let us clarify that we do not know whether the latter case is even possible, nor, if it is possible, whether the annealed probability of each limiting velocity must be $\frac12$, as one might be tempted to infer from our definition of $\vel(\mu)$. We could have defined $\vel(\mu)$ to be a convex combination of the two Dirac masses based on the annealed probability associated with each; however, this definition would make the statement of Corollary \ref{cor:continuity} weaker.} Note that if $\mu$ is a Dirac mass at a single vector $\xi\in\Sigma_{\mathcal{N}}^{\kappa,\eta,\theta}$, it follows from the strong law of large numbers that $\vel(\mu)$ is a Dirac mass at the vector $\left(\xi_i-\xi_{d+i}\right)_{i=1}^d\in B_1^1(0)$.

\begin{cor}
    If $\mu$ is a Dirac mass at a vector $\xi\in\Sigma_{\mathcal{N}}^{\kappa,\eta,\theta}$ with all positive entries, then $\vel$ is continuous at $\mu$. 
\end{cor}\label{cor:continuity}

\subsection{Main proof ideas}\label{subsec:ProofIdeas}

Let us describe the main ideas of the proofs of both Theorem \ref{thm:RecoveringKalikow+Bounds} and Theorem \ref{thm:MainThm}.
We define functions for a (finite or infinite) walk ${\bf X}=(X_n)_{n=0}^{\alpha}$, $\alpha\in\N\cup\{\infty\}$. For a vertex $x\in\Z^d$, define the \textit{hitting time} and \textit{positive hitting time} respectively as
\begin{equation}\notag
H_x({\bf X}):=\inf\{n\geq0:X_n=x\}\qquad\text{and}\qquad
\tilde{H}_x({\bf X}):=\inf\{n>0:X_n=x\}.
\end{equation}
and define the \textit{local time}
\begin{equation}\notag
N_x({\bf X}):=\#\{n\geq0:X_n=x\}.
\end{equation}
For a $T\in\N$, define the \textit{truncated local time}
\begin{equation}\notag
N_x^T({\bf X}):=\#\{0\leq n\leq T:X_n=x\}.
\end{equation}
As is typical with hitting times and other functions of a walk, we suppress the argument in the above functions when the walk referred to is clear from the context. The function $x\mapsto\Ea^0[N_x]$ is the \textit{annealed Green's function}, and a function  $x\mapsto\Ea^0[N_x^T]$ is a \textit{truncated annealed Green's function}. 

The method of our proof is to examine the truncated annealed \textit{Red's function} and \textit{Blue's function} given respectively by
\begin{equation}\notag
\Ea^0\left[N_x^T\mathbbm{1}_{\{x\text{ is red}\}}\right]\quad\quad\quad\text{and}\quad\quad\quad\Ea^0\left[N_x^T\mathbbm{1}_{\{x\text{ is blue}\}}\right].
\end{equation}
The goal is to show that that the truncated Red's function as a multiple of the truncated Blue's function can be bounded arbitrarily close to 0, uniformly in $x$ and in $T$, by making $p$ large enough. Once we do this, we can sum over all $x\in\Z^d$ and deduce that for any time $T$, the walk typically spends a very high proportion of the time from 0 to $T$ at blue sites, so that bounds on the quenched drift in a direction $\vec{u}\in S^{d-1}$ at blue sites ought to become near-bounds on the limiting speed in direction $\vec{u}$.  

This ``typical'' behavior, however, is given in terms of expectations rather than an almost-sure statement, and so more work is required to extract almost-sure bounds on the limiting speed. The arguments of \cite{Sznitman&Zerner1999} yield an almost sure limiting velocity, but one which may be random, supported on two points, in the case where a 0-1 law for directional transience does not hold.
Because the bound on the truncated Red's function is uniform for all $x\in\Z^d$, we may sum over all $x$ on each side of the hyperplane $\{x\cdot\vec{u}=0\}$ and get the same bounds. This argument essentially shows that conditioning on the walk's being transient in one direction or the other does not affect the bound on the expected proportion of time it spends at red sites, and therefore does not affect the bounds on the expected displacement after a large number of steps. Since there is at most one possible limiting velocity on each side of the hyperplane, these bounds on conditional expectations can then be used to obtain almost-sure bounds on the limiting velocity. The result of this work is Proposition \ref{prop:Bound}, which is used in the proof of Theorems \ref{thm:RecoveringKalikow+Bounds} and \ref{thm:MainThm}, reducing each of them to a matter of proving bounds, uniform in $x\in\Z^d$, on the truncated annealed Red's functions in terms of the truncated annealed Blue's functions.

For the proof of Theorem \ref{thm:RecoveringKalikow+Bounds}, our comparison of the  Red's and Blue's functions comes from exact expressions of ratios of these functions that involve probabilities of hitting and then returning to a site $x$, conditioned on the first step away from $x$. The main calculation is brief, but we need to truncate the walk at a random (geometric) time $\tau$, rather than a fixed $T$, in order to make it work. This geometric truncation is inspired by the methods of \cite{Sabot2004}.

Proving Theorem \ref{thm:MainThm}
is more difficult, and for this we use a coupling technique. Red sites are typically surrounded by blue sites. To show that a walk does not spend too much time on such red sites, we can couple two walks in environments that are exactly the same except at site $x$, which is blue in one environment and red in the other. The two walks decouple whenever they hit $x$. Then the particle in the environment for which $x$ is blue freezes, while the particle for which $x$ is red will hit $x$ at most a geometric number of times before it uses the blue sites surrounding $x$ to navigate to the point adjacent to $x$ where the first particle is waiting. At this time, they recouple and continue as before. This allows us to bound the expected amount of time at $x$ in the environment where $x$ is red and surrounded by blue sites in terms of the expected amount of time at $x$ in the environment where $x$ is blue, irrespective of the rest of the environment. Then, by taking $p$ large enough, we can make the environments where $x$ is red rare enough that under the annealed measure, the walk is expected to spend much more time in blue sites than in red sites surrounded by blue sites.

Of course, not every red site is completely surrounded by blue sites, and if $x$ is in a cluster of red sites, the above argument does not work directly. But we can use a similar argument to compare how much time the walk spends at sites $x$ that are in red clusters of various sizes, showing that the amount of time spent in a cluster of size $m$ decreases exponentially fast in $m$. 

\subsection*{Acknowledgments}

The author thanks Jonathon Peterson for useful conversations, for reviewing multiple versions, and for suggesting the idea of following \cite{Sabot2004} in using a geometric truncation for the proof of Theorem \ref{thm:RecoveringKalikow+Bounds}. Some of the work on this paper was supported through NSF grant DMS-2153869. The author is grateful to Leonid Petrov for providing this support.

\section{Preliminaries}\label{sec:Prelim}

We describe a condition that will replace a traditional ellipticity assumption.

\begin{enumerate}[(C)]
    \item The probability $p$ and the measure $\mu_r$ are\footnote{Note that we do not include anything about $\mu_b$ in condition (C) because $\mu_b$ is already assumed to satisfy a uniform ellipticity assumption.} such that with $\Pm$-probability 1, the Markov chain induced by $\omega$ has only one infinite communicating class, and it is reachable from every site.
\end{enumerate}

This is clearly satisfied under the traditional ellipticity assumption where every vector in the support of $\mu_r$ has all positive components, but it is also satisfied if, for example, $p>p_c(d)$, where $p_c(d)$ is the critical site percolation threshold for $\Z^d$.

\begin{lem}\label{lem:condC}
    Condition (C) is satisfied if $p>p_c(d)$.
\end{lem}
\begin{proof}
    Assume $p>p_c(d)$. Then it is standard that there is a unique infinite cluster\footnote{Elsewhere in this paper, we use the word ``cluster'' in our own specialized sense, but in this argument, we use it in the traditional sense for percolation on $\Z^d$.} of blue sites (see, e.g., \cite[Theorem 2.4]{Hoggstrom&Jonasson2006}, which is given for bond percolation, but comments on pages 332 and 333 explain that the proof adapts with no additional complication to the case of site percolation). All sites in the infinite blue cluster belong to a single communicating class by the ellipticity assumption on blue sites.
    Since the origin has positive probability of belonging to the infinite blue cluster, the ergodic theorem implies that with probability 1, there are sites in the infinite blue cluster on both sides of the origin on the discrete line $x+\Z_{\eta}$ for any $x\in\Z^d$. Thus, from every $x\in\Z^d$ it is possible to reach the infinite blue cluster by taking steps in direction $e_{\eta}$, and every $x\in\Z^d$ can be reached in the infinite blue cluster by taking steps in direction $e_{\eta}$. Thus, all sites almost surely belong to a single communicating class, which is stronger than condition (C).
\end{proof}

Under Condition (C), we have the following law of large numbers.
 
    \begin{thm}\label{claim:449}
        Assume Condition (C). Then $\frac{X_n\cdot\vec{u}}{n}$ has a $\Pa^0$--almost sure limit. The limit, \textit{a priori}, could be random, but one of the following situations holds:
        \begin{enumerate}
            \item $\Pa^0(A_{\vec{u}}\cup A_{-\vec{u}})=0$ and with probability 1, $\lim_{n\to\infty}\frac{X_n\cdot\vec{u}}{n}=0$;
            \item $\Pa^0(A_{\vec{u}})=1$ and $\lim_{n\to\infty}\frac{X_n\cdot\vec{u}}{n}$ is deterministic;
            \item $\Pa^0(A_{-\vec{u}})=1$ and $\lim_{n\to\infty}\frac{X_n\cdot\vec{u}}{n}$ is deterministic;
            \item $\Pa^0(A_{\vec{u}}\cup A_{-\vec{u}})=1$, $0<\Pa^0(A_{\vec{u}})<1$, and $\lim_{n\to\infty}\frac{X_n\cdot\vec{u}}{n}$ is supported on at most two values, $L_1\leq0$ and $L_2\geq0$, such with probability 1, if $\lim_{n\to\infty}X_n\cdot\vec{u}=-\infty$, then $\lim_{n\to\infty}\frac{X_n\cdot\vec{u}}{n}=L_1$, and if $\lim_{n\to\infty}X_n\cdot\vec{u}=\infty$, then $\lim_{n\to\infty}\frac{X_n\cdot\vec{u}}{n}=L_2$.
        \end{enumerate}
    \end{thm}
    This is the main theorem of \cite{Sznitman&Zerner1999} and \cite{Zerner2002}. See in particular the remark in \cite{Zerner2002} after the statement of Theorem A, noting that the uniform ellipticity assumption assumed in \cite{Sznitman&Zerner1999} is unnecessary, and is not used in the proof given in \cite{Zeitouni2004}. Indeed, the proofs do not even use the weaker ellipticity assumption that $\omega^0(e_i)>0$, $\Pm$--a.s., for all $1\leq i\leq 2d$, except insofar as they quote Lemma 4 and Proposition 3 of \cite{Zerner&Merkl2001}. But these results are proven in \cite{Slonim2021b} as Lemma 3.5 and Theorem A.1, respectively, under condition (C) only. We may therefore take Theorem \ref{claim:449} as proven.  

\section{Using Blue's and Red's functions to bound limiting velocity}\label{sec:Criterion}

In this section, we give a criterion for directional transience and ballisticity in terms of truncated Blue's and Red's functions. The criterion will be used in the proofs of both Theorem \ref{thm:MainThm} and Theorem \ref{thm:RecoveringKalikow+Bounds}. Because the latter proof requires truncating at a random time, we give the criterion in terms of a geometric random variable $\tau$ rather than a fixed $T$. 

\begin{prop}\label{prop:Bound}
    Assume condition (C). Fix a unit vector $\vec{u}\in\R^d$. Suppose $v_1^b\geq v_1^r$ are such that for every $\xi$ in the support of $\mu_b$,  $\sum_{i=1}^{2d}\xi(e_i)e_i\cdot\vec{u}\geq v_1^b$, and for every $\xi$ in the support of $\mu_r$,
    $\sum_{i=1}^{2d}\xi(e_i)e_i\cdot\vec{u}\geq v_1^r$.
        Let $\tau$ be a geometric random variable, independent of everything else, with expectation $\frac{1}{\rho}$. If there exists $\alpha>0$ such that
    \begin{equation}\label{eqn:395}    
    \Ea^0[N_x^{\tau}\mathbbm{1}_{\{x\text{ is blue}\}}]\geq\alpha\Ea^0[N_x^{\tau}\mathbbm{1}_{\{x\text{ is red}\}}]
    \end{equation}
    for all $x\in\Z^d$ and for all sufficiently small $\rho$, then with $\Pa^0$-probability 1,
    \begin{equation}
        \lim_{n\to\infty}\frac{X_n\cdot\vec{u}}{n}\geq \frac{\alpha}{\alpha+1}v_1^b+\frac{1}{\alpha+1}v_1^r.
    \end{equation}
\end{prop}

The proof will use the following notation: For a unit vector $\vec{u}\in\R^d$, let $A_{\vec{u}}$ be the event that the walk is transient in direction $\vec{u}$; i.e., $\lim_{n\to\infty}X_n\cdot\vec{u}=\infty$.

\begin{proof} If necessary, expand the probability space to accommodate $\tau$. 
We divide the proof into a series of claims. 

    \begin{claim}\label{claim:403}
    \begin{align}\label{eqn:440}
        \frac{\Ea^0[X_{\tau}\cdot\vec{u}]}{\Ea^0[\tau]}
            \geq \frac{\alpha}{\alpha+1}v_1^b+\frac{1}{\alpha+1}v_1^r.
    \end{align}
    \end{claim}
    To prove this claim, we begin with
    \begin{align*}\notag
        X_{\tau}\cdot\vec{u}&=\sum_{j=0}^{\tau-1}(X_{j+1}-X_j)\cdot\vec{u}
        \\
        &=\sum_{j=0}^{\infty}\mathbbm{1}_{\{j<\tau\}}\mathbbm{1}_{\{\text{$x$ is blue}\}}(X_{j+1}-X_j)\cdot\vec{u}
        +
        \sum_{j=0}^{\infty}\mathbbm{1}_{\{j<\tau\}}\mathbbm{1}_{\{\text{$x$ is red}\}}(X_{j+1}-X_j)\cdot\vec{u}.
    \end{align*}
    Taking expectations, we get 
        \begin{align}\notag
            \Ea^0[X_{\tau}\cdot\vec{u}]&=\sum_{j=0}^{\infty}\Pa^0(j<\tau)\Pa^0(X_j\text{ is blue})\Ea^0[(X_{j+1}-X_j)\cdot\vec{u}\givenn X_j\text{ is blue}]
            \\\notag
            &\qquad\qquad
            +\sum_{j=0}^{\infty}\Pa^0(j<\tau)\Pa^0(X_j\text{ is red})\Ea^0[(X_{j+1}-X_j)\cdot\vec{u}\givenn X_j\text{ is red}]
            \\\notag
            &\geq\sum_{j=0}^{\infty}\Pa^0(j<\tau)\Pa^0(X_j\text{ is blue})v_1^b
            +\sum_{j=0}^{\infty}\Pa^0(j<\tau)\Pa^0(X_j\text{ is red})v_1^r
            \\\label{eqn:528}
            &=v_1^b\sum_{x\in\Z^d}\Ea^0[N_x^{\tau}\mathbbm{1}_{\{x\text{ is blue}\}}]
            +v_1^r\sum_{x\in\Z^d}\Ea^0[N_x^{\tau}\mathbbm{1}_{\{x\text{ is red}\}}]
        \end{align}
        
        Now if $a,b>0$, $c>d$, and $a>\alpha b$, then $\frac{ac+bd}{a+b}\geq\frac{\alpha}{\alpha+1}c+\frac{1}{\alpha+1}d$. 
        Letting $a=\sum_{x\in\Z^d}\Ea^0\left[N_x^{\tau}\mathbbm{1}_{\{\text{$x$ is blue}\}}\right]$, $b=\sum_{x\in\Z^d}\Ea^0\left[N_x^{\tau}\mathbbm{1}_{\{\text{$x$ is red}\}}\right]$, $c=v_1^b$, and $d=v_1^r$, we apply this fact along with \eqref{eqn:395}, \eqref{eqn:528}, and the fact that 
        \begin{equation}\notag
        \sum_{x\in\Z^d}\Ea^0[N_x^{\tau}\mathbbm{1}_{\{x\text{ is blue}\}}]+\sum_{x\in\Z^d}\Ea^0[N_x^{\tau}\mathbbm{1}_{\{x\text{ is red}\}}]=\Ea^0[\tau]
        \end{equation}
        to get \eqref{eqn:440}, proving Claim \ref{claim:403}. Now for convenience, we define
        $$
        v':=\frac{\alpha}{\alpha+1}v_1^b+\frac{1}{\alpha+1}v_1^r
        $$
        to be the right side of \eqref{eqn:440}.
        The next claim turns information about $X_{\tau}$ into information about $X_n$.
        \begin{claim}\label{claim:436} 
        \begin{equation}\notag
            \limsup_{n\to\infty}\frac{\Ea^0[X_n\cdot\vec{u}]}{n}>v'.
        \end{equation} 
        \end{claim}
    We prove this by contradiction. For if there is a $c<v'$ such that eventually $\frac{\Ea^0[X_n\cdot\vec{u}]}{n}\leq c$ (say for $n\geq L$), then taking 
    \begin{align}\notag
        \Ea^0[X_{\tau}\cdot \vec{u}]&=\sum_{j=0}^{\infty}\Pa^0(\tau=j)\Ea^0[X_j\cdot\vec{u}]
        \\\notag
        &\leq\sum_{j=L}^{\infty}\Pa^0(\tau=j)cj+\Pa^0(\tau<L)\max_{0\leq j\leq L}\Ea^0[X_j\cdot\vec{u}]
        \\\notag
        &\leq c\Ea^0[\tau]+c',
    \end{align}
where $c'=\Pa^0(\tau<L)\max_{0\leq j\leq L}\Ea^0[X_j\cdot\vec{u}]$, and by taking $\tau$ to have small enough parameter, we may make $c'$ as small as we please. In particular, since $c<v'$, we may take $c'$ small enough that $c+c'<v'$, and get $\Ea^0[X_{\tau}\cdot\vec{u}]\leq (c+c')\Ea^0[\tau]<v'\Ea^0[\tau]$, contradicting Claim \ref{claim:403}. This proves Claim \ref{claim:436}.

    Now if $\lim_{n\to\infty}\frac{X_n\cdot\vec{u}}{n}$ is deterministic, then Claim \ref{claim:436} finishes the proof. We may therefore assume we are in the fourth situation from Theorem \ref{claim:449}: namely, that $\Pa^0(A_{\vec{u}}\cup A_{-\vec{u}})=1$, $0<\Pa^0(A_{\vec{u}})<1$, and $\lim_{n\to\infty}\frac{X_n\cdot\vec{u}}{n}$ is supported on at most two values, $L_1\leq0$ and $L_2\geq0$, such with probability 1, if $\lim_{n\to\infty}X_n\cdot\vec{u}=-\infty$, then $\lim_{n\to\infty}\frac{X_n\cdot\vec{u}}{n}=L_1$, and if $\lim_{n\to\infty}X_n\cdot\vec{u}=\infty$, then $\lim_{n\to\infty}\frac{X_n\cdot\vec{u}}{n}=L_2$.

\begin{claim}\label{claim:598}
    \begin{equation}
        \lim_{\rho\to0}\frac{\Ea^0\left[X_{\tau}\cdot\vec{u}\mathbbm{1}_{A_{-\vec{u}}}\right]}{\Ea^0\left[\tau\mathbbm{1}_{A_{-\vec{u}}}\right]}=L_1
    \end{equation}
\end{claim}

To prove this claim, note that from Theorem \ref{claim:449}, we have that on the event $A_{-\vec{u}}$, $\lim_{n\to\infty}\frac{X_n}{n}=L_1$. Let $\delta>0$. We can choose $N$ such that the probability, conditioned on $A_{-\vec{u}}$, that there exists $n\geq N$ with $\left|\frac{X_n\cdot\vec{u}}{n}-L_1\right|\geq\frac{\delta}{5}$ is less than $\frac{\delta}{5}$. Then we may choose $\rho$ small enough that the probability that $\tau<N$ is less than $\frac{\delta}{5}$. For such a $\rho$, it is straightforward to check that we have
\begin{equation*}
\left|\frac{\Ea^0\left[X_{\tau}\cdot\vec{u}\mathbbm{1}_{A_{-\vec{u}}}\right]}{\Ea^0\left[\tau\mathbbm{1}_{A_{-\vec{u}}}\right]}-L_1\right|<\delta.
\end{equation*}

This proves Claim \ref{claim:598}.

\begin{claim}\label{claim:604}
    \begin{equation}\notag
        \lim_{\rho\to0}\frac{\Ea^0\left[\sum_{j=1}^{\tau}\mathbbm{1}_{\{X_{j-1}\cdot\vec{u}\leq0\}}\mathbbm{1}_{A_{\vec{u}}}\right]}{\Ea^0[\tau]}=0
        \qquad\text{and}\qquad
        \lim_{\rho\to0}\frac{\Ea^0\left[\sum_{j=1}^{\tau}\mathbbm{1}_{\{X_{j-1}\cdot\vec{u}>0\}}\mathbbm{1}_{A_{-\vec{u}}}\right]}{\Ea^0[\tau]}=0
    \end{equation}
\end{claim}

We prove that the first fraction approaches 0; the second follows by an analogous argument. Let $\delta>0$. On the event $A_{\vec{u}}$, there exists some $K$ such that for $n\geq K$, $X_n\cdot\vec{u}>0$. Therefore, we may choose $K$ such that the probability that $X_n\cdot\vec{u}\leq 0$ for some $n\geq K$ and $A_{\vec{u}}$ holds is less than $\frac{\delta}{2}$. Let $A_{\vec{u}}^{K-\text{bad}}$ be this unlikely event. Now let $\rho<\frac{\delta}{2K}$, so that $K\leq \frac{\delta}{2}\Ea^0[\tau]$. 
Then
\begin{align*}
    \sum_{j=1}^{\tau}\mathbbm{1}_{\{X_{j-1}\cdot\vec{u}\leq0\}}\mathbbm{1}_{A_{\vec{u}}}
    &\leq
    K+
    \tau\mathbbm{1}_{A_{\vec{u}}^{K-\text{bad}}}
    \\
    &\leq 
    \frac{\delta}{2}\Ea^0[\tau]
    +
    \Ea^0[\tau]\Pa^0(A_{\vec{u}}^{K-\text{bad}})
    \\
    &<\delta\Ea^0[\tau].
\end{align*}
This proves the first half of the statement of Claim \ref{claim:604}, and the second half is proved in an analogous way.

\begin{claim}\label{claim:571}
\begin{equation}\label{eqn:612}
    \liminf_{\rho\to0}\frac{\Ea^0\left[X_{\tau}\cdot\vec{u}\mathbbm{1}_{A_{-\vec{u}}}\right]}{\Ea^0\left[\tau\mathbbm{1}_{A_{-\vec{u}}}\right]}\geq v'.
\end{equation}
\end{claim}

The numerator of the left hand side of \eqref{eqn:612} is
\begin{align*}
    \Ea^0\left[X_{\tau}\cdot\vec{u}\mathbbm{1}_{A_{-\vec{u}}}\right]
    &=\Ea^0\left[\sum_{j=1}^{\tau}(X_{j}-X_{j-1})\cdot\vec{u}\mathbbm{1}_{\{X_{j-1}\cdot\vec{u}\leq0\}}\mathbbm{1}_{A_{-\vec{u}}}\right]+\Ea^0\left[\sum_{j=1}^{\tau}(X_{j}-X_{j-1})\cdot\vec{u}\mathbbm{1}_{\{X_{j-1}\cdot\vec{u}>0\}}\mathbbm{1}_{A_{-\vec{u}}}\right]
    \\
    &\geq\Ea^0\left[\sum_{j=1}^{\tau}(X_{j}-X_{j-1})\cdot\vec{u}\mathbbm{1}_{\{X_{j-1}\cdot\vec{u}\leq0\}}\right]
    \\
    &\qquad\qquad
    -\Ea^0\left[\sum_{j=1}^{\tau}\mathbbm{1}_{\{X_{j-1}\cdot\vec{u}\leq0\}}\mathbbm{1}_{A_{\vec{u}}}\right]-\Ea^0\left[\sum_{j=1}^{\tau}\mathbbm{1}_{\{X_{j-1}\cdot\vec{u}>0\}}\mathbbm{1}_{A_{-\vec{u}}}\right]
\end{align*}
    By Claim \ref{claim:604}, the last two terms vanish as a fraction of $\Ea^0[\tau]$ as $\rho\to\infty$, and thus also vanish as a fraction of $\Ea^0[\tau\mathbbm{1}_{A_{-\vec{u}}}]=\Ea^0[\tau]\Pa^0(A_{-\vec{u}})$, since by assumption $\Pa^0(A_{-\vec{u}})>0$.

    We also have
    \begin{align*}
    \Ea^0\left[\tau\mathbbm{1}_{A_{-\vec{u}}}\right]
    &\leq\Ea^0\left[\sum_{j=1}^{\tau}\mathbbm{1}_{\{X_{j-1}\cdot\vec{u}\leq0\}}\right]
    +\Ea^0\left[\sum_{j=1}^{\tau}\mathbbm{1}_{\{X_{j-1}\cdot\vec{u}>0\}}\mathbbm{1}_{A_{-\vec{u}}}\right],
\end{align*}
and the second term vanishes as a fraction of $\Ea^0[\tau\mathbbm{1}_{A_{-\vec{u}}}]$.
Therefore,
\begin{align}\label{eqn:663}
    \liminf_{\rho\to0}\frac{\Ea^0\left[X_{\tau}\cdot\vec{u}\mathbbm{1}_{A_{-\vec{u}}}\right]}{\Ea^0\left[\tau\mathbbm{1}_{A_{-\vec{u}}}\right]}
    &\geq
    \liminf_{\rho\to0}\frac{\Ea^0\left[\sum_{j=1}^{\tau}(X_{j}-X_{j-1})\cdot\vec{u}\mathbbm{1}_{\{X_{j-1}\cdot\vec{u}\leq0\}}\right]}{\Ea^0\left[\sum_{j=1}^{\tau}\mathbbm{1}_{\{X_{j-1}\cdot\vec{u}\leq0\}}\right]}.
\end{align}
Let us examine the numerator on the right in the above inequality.
\begin{align}\notag
\Ea^0\left[\sum_{j=1}^{\tau}(X_j-X_{j-1})\cdot\vec{u}\mathbbm{1}_{\{X_{n-1}\cdot\vec{u}\leq0\}}
\right]
&=
\sum_{x\cdot\vec{u}\leq0}\Ea^0\left[\sum_{j=1}^{\tau}(X_j-X_{j-1})\mathbbm{1}_{\{X_{n-1}=x\}}\right]
\\\notag
&=
\sum_{x\cdot\vec{u}\leq0}\Ea^0\left[\sum_{j=1}^{\tau}(X_j-X_{j-1})\mathbbm{1}_{\{X_{n-1}=x\}}\mathbbm{1}_{\{x\text{ is blue}\}}\right]
\\\notag
&\qquad+
\sum_{x\cdot\vec{u}\leq0}\Ea^0\left[\sum_{j=1}^{\tau}(X_j-X_{j-1})\mathbbm{1}_{\{X_{n-1}=x\}}\mathbbm{1}_{\{x\text{ is red}\}}\right]
\\\notag
&\geq
\sum_{x\cdot\vec{u}\leq0}v_1^b\Ea^0\left[\sum_{j=1}^{\tau}\mathbbm{1}_{\{X_{n-1}=x\}}\mathbbm{1}_{\{x\text{ is blue}\}}\right]
\\\notag
&\qquad+
\sum_{x\cdot\vec{u}\leq0}v_1^r\Ea^0\left[\sum_{j=1}^{\tau}\mathbbm{1}_{\{X_{n-1}=x\}}\mathbbm{1}_{\{x\text{ is red}\}}\right]
\\\label{eqn:689}
&=v_1^b\sum_{x\cdot\vec{u}\leq0}\Ea^0[N_x^{\tau-1}\mathbbm{1}_{\{x\text{ is blue}\}}]+v_1^r\sum_{x\cdot\vec{u}\leq0}\Ea^0[N_x^{\tau-1}\mathbbm{1}_{\{x\text{ is red}\}}]
\end{align}
The denominator is
\begin{equation}\label{eqn:693}
    {\Ea^0\left[\sum_{j=1}^{\tau}\mathbbm{1}_{\{X_{j-1}\cdot\vec{u}\leq0\}}\right]}
    =
    \sum_{x\cdot\vec{u}\leq0}\Ea^0[N_x^{\tau-1}\mathbbm{1}_{\{x\text{ is blue}\}}]+\sum_{x\cdot\vec{u}\leq0}\Ea^0[N_x^{\tau-1}\mathbbm{1}_{\{x\text{ is red}\}}]
\end{equation}

        Again we use the fact that  if $a,b>0$, $c>d$, and $a>\alpha b$, then $\frac{ac+bd}{a+b}\geq\frac{\alpha}{\alpha+1}c+\frac{1}{\alpha+1}d$. Letting $a=\sum_{x\cdot\vec{u}\leq0}\Ea^0\left[N_x^{\tau}\mathbbm{1}_{\{\text{$x$ is blue}\}}\right]$, $b=\sum_{x\cdot\vec{u}\leq0}\Ea^0\left[N_x^{\tau}\mathbbm{1}_{\{\text{$x$ is red}\}}\right]$,  $c=v_1^b$, and $d=v_1^r$ we apply this fact along with \eqref{eqn:395} \eqref{eqn:689}, and \eqref{eqn:693} to get
\begin{equation}
    \frac{\Ea^0\left[\sum_{j=1}^{\tau}(X_{j}-X_{j-1})\cdot\vec{u}\mathbbm{1}_{\{X_{j-1}\cdot\vec{u}\leq0\}}\right]}{\Ea^0\left[\sum_{j=1}^{\tau}\mathbbm{1}_{\{X_{j-1}\cdot\vec{u}\leq0\}}\right]}>v'.
\end{equation}
Combining this with \eqref{eqn:663} finishes the proof of Claim \ref{claim:571}. 

Combining Claim \ref{claim:571} with Claim \ref{claim:598} proves that $L_1\geq v'$, which finishes the proof of Proposition \ref{prop:Bound}.
\end{proof}

We now state a simple corollary to Proposition \ref{prop:Bound}.
 
\begin{cor}\label{cor:DirectionalTransience}
   Assume condition (C). Fix a unit vector $\vec{u}\in\R^d$. Suppose $v_1^b,v_2^b$ are such that for every $\xi$ in the support of $\mu_b$, $\sum_{i=1}^{2d}\xi(e_i)e_i\cdot\vec{u}\in[v_1,v_2]$. Let $\varepsilon>0$. If 
    \begin{equation}\label{eqn:395c}    
    \Ea^0[N_x^{T}\mathbbm{1}_{\{x\text{ is blue}\}}]\geq\frac{2-\varepsilon}{\varepsilon}\Ea^0[N_x^{T}\mathbbm{1}_{\{x\text{ is red}\}}]
    \end{equation}
    for all $x\in\Z^2$ and for all sufficiently large $T$, then with $\Pa^0$-probability 1,
    \begin{equation}\label{eqn:635}
        v_1^b-\varepsilon\leq\lim_{n\to\infty}\frac{X_n\cdot\vec{u}}{n}\leq v_2^b+\varepsilon.
    \end{equation}
\end{cor}

\begin{proof}
    Let $\varepsilon>0$. By the nearest-neighbor assumption, we necessarily have $\frac{X_n\cdot\vec{u}}{n}\in[-1,1]$. We therefore assume $v_1,v_2\in[-1,1]$ and $\varepsilon<2$, since the statement is trivial otherwise. We may define $v_1^r=-1$ so that $v_1^b\geq v_1^r$.  
    If the assumptions here are met, then by conditioning on the value of $\tau$, the assumptions of Proposition \ref{prop:Bound} are met with $\alpha=\frac{2-\varepsilon}{\varepsilon}$. We conclude that with  $\Pa^0$-probability 1,
    \begin{align*}
        \lim_{n\to\infty}\frac{X_n\cdot\vec{u}}{n}&\geq \frac{\frac{2-\varepsilon}{\varepsilon}}{\frac{2-\varepsilon}{\varepsilon}+1}v_1^b+\frac{1}{\frac{2-\varepsilon}{\varepsilon}+1}v_1^r
        \\
        &=\left(1-\frac{\varepsilon}{2}\right)v_1^b-\frac{\varepsilon}{2}
        \\
        &=v_1-\varepsilon.
    \end{align*}
    This gives us the first inequality in \eqref{eqn:635}. By applying it in direction $-\vec{u}$, with $-v_2$ playing the role of $v_1$, we get the second inequality.
\end{proof}

\section{Recovering Kalikow's theorem with quantitative ballisticity bounds}\label{sec:Kalikow}

In this section, we use Proposition \ref{prop:Bound} to prove Theorem \ref{thm:RecoveringKalikow+Bounds}. In order to recover the full strength of Kalikow's result, we will need to truncate the walk at a geometric time instead of a deterministic time. To do this, we must expand our probability space. We previously have defined $P_{\omega}^x$ to be a measure on $(\Z^d)^{\N}$. Now we let $P_{\omega}^x$ be the measure on $(\Z^d)^{\N}\times\N$ giving the law of a random walk in $\omega$ started at $x$ and an independent geometric random variable $\tau$ with expectation $\frac{1}{\rho}$. We do not change the measure $\Pm$, and as before, we let $\Pa^x=\Pm\times P_{\omega}^x$---formally defined by \eqref{eqn:362}---be the annealed measure starting at $x$.

\begin{lem}\label{lem:1008}
    In Kalikow's two-vertex model, for all $x\in\Z^2$,
\begin{equation}\notag
    \Ea^0\left[
    N_x^{\tau}\one{$x$ is red}
    \right]
    \leq \frac{(1-p)M}p
    \Ea^0\left[
    N_x^{\tau}\one{$x$ is blue}
    \right],
\end{equation}
where $M$ is as  defined in \eqref{eqn:220}. 
\end{lem}

\begin{proof}
Let  $\gamma=1-\rho$ be the one-step survival probability. Then 
\begin{align}\notag
     \notag\frac{\Ea^0[N_x^{\tau}\one{$x$ is red}]}{\Ea^0[N_x^{\tau}\one{$x$ is blue}]}
     &=
     \frac{\Em[E_{\omega}^0[N_x^{\tau}]\mathbbm{1}_{x\text{ is red}}]}{\Em[E_{\omega}^0[N_x^{\tau}]\mathbbm{1}_{x\text{ is blue}}]}
     \\\notag
     &=
     \frac{\Em\left[(P_{\omega}^0(H_x<\tau)E_{\omega}^0[N_x^{\tau}\given H_x<\tau]+P_{\omega}^0(H_x=\tau))\one{$x$ is red}\right]}{\Em\left[(P_{\omega}^0(H_x<\tau)E_{\omega}^0[N_x^{\tau}\given H_x<\tau]+P_{\omega}^0(H_x=\tau))\one{$x$ is blue}\right]}
     \\\notag
     &=
     \frac{\Em\left[(P_{\omega}^0(H_x<\tau)E_{\omega}^x[N_x^{\tau}]+P_{\omega}^0(H_x=\tau))\one{$x$ is red}\right]}{\Em\left[(P_{\omega}^0(H_x<\tau)E_{\omega}^x[N_x^{\tau}]+P_{\omega}^0(H_x=\tau))\one{$x$ is blue}\right]}
     \\\notag
     &=
     \frac{\Em\left[(P_{\omega}^0(H_x<\tau)E_{\omega}^x[N_x^{\tau}]\one{$x$ is red}\right]+(1-p)\Pa^0(H_x=\tau)}{\Em\left[(P_{\omega}^0(H_x<\tau)E_{\omega}^x[N_x^{\tau}]\one{$x$ is blue}\right]+p\Pa^0(H_x=\tau)}.
     \end{align}
To get the third equality, we used memorylessness of geometric random variables. The fourth is obtained by distributing the indicator and noticing that the walk up until $H_x$ is independent of the color of $x$.
Now for fixed $\omega$, $E_{\omega}^x[N_x^{\tau}]$ is the expectation of a geometric random variable with parameter $P_{\omega}^x(\tau<\tilde{H}_x)$. Using this fact and letting $a=\Pa^0(H_x=\tau)$. , we get
     \begin{align}\notag
     \notag\frac{\Ea^0[N_x^{\tau}\one{$x$ is red}]}{\Ea^0[N_x^{\tau}\one{$x$ is blue}]}
     &=
     \frac{\Em\left[\frac{P_{\omega}^0(H_x<\tau)}{P_{\omega}^x(\tau<\tilde{H}_x)}\one{$x$ is red}\right]+(1-p)a}{\Em\left[\frac{P_{\omega}^0(H_x<\tau)}{P_{\omega}^x(\tau<\tilde{H}_x)}\one{$x$ is blue}\right]+pa}
     \\\notag
     &=
     \frac{\Em\left[\frac{P_{\omega}^0(H_x<\tau)}{\sum_{i=1}^4\gamma r_iP_{\omega}^{x+e_i}(\tau<H_x)}\right](1-p)+(1-p)a}
     {\Em\left[\frac{P_{\omega}^0(H_x<\tau)}{\sum_{i=1}^4\gamma b_iP_{\omega}^{x+e_i}(\tau<H_x)}\right]p+pa}
     \\\label{eqn:237}
     &=\frac{1-p}{p}\frac{\Em[Y+a]}{\Em[Z+a]},
\end{align}
where $Y$ and $Z$ are the fractions inside the expectations in the second to last line. Now with probability 1 we have
\begin{align}\notag
\frac{Y}{Z}&=\frac{\sum_{i=1}^4\gamma b_iP^{x+e_i}(\tau<H_x)}{\sum_{i=1}^4\gamma r_iP^{x+e_i}(\tau<H_x)}
\end{align}
Let $q_i=P_{\omega}^{x+e_i}(\tau<H_x)$. Then almost surely
\begin{align}\notag
    \frac{Y}{Z}
    &=\frac{\sum_{i=1}^4\gamma b_iq_i}{\sum_{i=1}^4\gamma r_iq_i},
    \\\notag
    &\leq\max_{i=1,2,3,4}\frac{b_i}{q_i}=M,
\end{align}
using the fact that a ratio of weighted sums is no more than the maximum ratio of individual terms.
Therefore, $\Em[Y]\leq M\Em[Z]$. Because $M\geq1$, this implies that $\Em[Y+a]\leq M\Em[Z+a]$, so by \eqref{eqn:237}, $\frac{\Ea^0[N_x^{\tau}\one{$x$ is red}]}{\Ea^0[N_x^{\tau}\one{$x$ is blue}]}\leq \frac{1-p}{p}M$, which is what we wanted to show.
\end{proof}

\begin{proof}[Proof of Theorem \ref{thm:RecoveringKalikow+Bounds}]
Note that in the two-vertex model, we can use $b_1-b_3$ and $r_3-r_1$ as the $v_1^b$ and $v_r$ in Proposition \ref{prop:Bound}.
If we define 
\[
\alpha:=\frac{p}{(1-p)M},
\]
then $\alpha>0$, and Lemma \ref{lem:1008} gives us
\begin{align}\notag
    \Ea^0\left[
    N_x^{\tau}\one{$x$ is blue}
    \right]
    &\geq \frac{p}{(1-p)M}
    \Ea^0\left[
    N_x^{\tau}\one{$x$ is red}
    \right]
    \\\notag
    &=\alpha\Ea^0\left[
    N_x^{\tau}\one{$x$ is red}
    \right]
\end{align}
Now the assumptions of Theorem \ref{thm:RecoveringKalikow+Bounds} imply ellipticity, which implies condition (C). Therefore, Proposition \ref{prop:Bound} tells us that $\Pa^0$--a.s.,
    \begin{align*}
        \lim_{n\to\infty}\frac{X_n\cdot e_1}{n}&\geq
        \frac{\alpha}{\alpha+1}(b_1-b_3)+\frac{1}{\alpha+1}(r_3-r_1)
        \\
        &=
    \frac{\frac{p}{(1-p)M}}{\frac{p}{(1-p)M}+1}(b_1-b_3)+\frac{1}{\frac{p}{(1-p)M}+1}(r_3-r_1)
    \\
    &=\frac{p(b_1-b_3)-(1-p)M(r_3-r_1)}{p+(1-p)M}.
    \end{align*}
Now \eqref{eqn:220} implies that this lower bound is positive, so that the walk is almost surely transient and ballistic in direction $e_1$. By Theorem \ref{claim:449}, this means the limiting velocity is deterministic.

We now seek to prove the upper bound. For $k\in\N$, let $\text{New}_k$ be the event that $X_k\notin\{X_0,X_1,\ldots,X_{k-1}\}$; that is, that the walk reaches a ``new'' site at time $k$. By the nearest-neighbor assumption, we have
    \begin{equation}
        X_n\cdot e_1\leq \sum_{k=1}^{n}\mathbbm{1}_{\text{New}_k}\leq \sum_{k=0}^{n-1}\mathbbm{1}_{\text{New}_k},
    \end{equation}
    so
    \begin{equation}\label{eqn:1369}
    \lim_{n\to\infty}\frac{\Ea^0[X_n\cdot e_1]}n\leq \liminf_{n\to\infty}\frac1n\sum_{k=0}^{n-1}\Pa^0(\text{New}_k)
    \end{equation}
On the other hand, 
\begin{align}
    X_n\cdot e_1&=
    \sum_{k=0}^{n-1}(X_{k+1}-X_k)\cdot e_1\mathbbm{1}_{\text{New}_k}
    +
    \sum_{k=0}^{n-1}(X_{k+1}-X_k)\cdot e_1\mathbbm{1}_{\text{New}_k^c},
\end{align}
    so
\begin{align*}
    \lim_{n\to\infty}\frac{\Ea^0[X_n\cdot e_1]}n&\leq
    \liminf_{n\to\infty}\frac1n\sum_{k=0}^{n-1}\Pa^0(\text{New}_k)\Ea^0[(X_{k+1}-X_k)\cdot e_1|\text{New}_k]
    \\
    &\qquad\qquad
    +\liminf_{n\to\infty}\frac1n\sum_{k=0}^{n-1}\Pa^0(\text{New}_k^c)\Ea^0[(X_{k+1}-X_k)\cdot e_1|\text{New}_k^c]
    \\
    &\leq
    \liminf_{n\to\infty}\frac1n\sum_{k=0}^{n-1}\Pa^0(\text{New}_k)[p(b_1-b_3)-(1-p)(r_3-r_1)]
    \\
    &\qquad\qquad
    +\limsup_{n\to\infty}\frac1n\sum_{k=0}^{n-1}\Pa^0(\text{New}_k^c)(b_1-b_3).
    \\
    &=
    \left(\liminf_{n\to\infty}\frac1n\sum_{k=0}^{n-1}\Pa^0(\text{New}_k)\right)[p(b_1-b_3)-(1-p)(r_3-r_1)]
    \\
    &\qquad\qquad
    +\left(1-\liminf_{n\to\infty}\frac1n\sum_{k=0}^{n-1}\Pa^0(\text{New}_k)\right)(b_1-b_3).
\end{align*}
The second inequality comes from the fact that conditioning on being at a fresh site at time $k$ (and conditioning on no information about the future of the walk), the expected displacement from $X_k$ to $X_{k+1}$ is simply the annealed drift, and also that conditioning on any information that does not include information about the future of the walk, the expected displacement from  $X_k$ to $X_{k+1}$ is no more than the maximum possible quenched drift. Now the last line (which we have broken into two lines) is a convex combination of $[p(b_1-b_3)-(1-p)(r_3-r_1)]$ and $(b_1-b_3)$, the latter being the larger quantity. Increasing the weight given to the larger quantity increases the value of the convex combination, so using \eqref{eqn:1369} lets us continue:

\begin{align*}
    \lim_{n\to\infty}\frac{\Ea^0[X_n\cdot e_1]}n&\leq
    \left(\lim_{n\to\infty}\frac{\Ea^0[X_n\cdot e_1]}n\right)[p(b_1-b_3)-(1-p)(r_3-r_1)]
    \\
    &\qquad\qquad
    +\left(1-\lim_{n\to\infty}\frac{\Ea^0[X_n\cdot e_1]}n\right)(b_1-b_3).
\end{align*}

Solving, we get
$$
\lim_{n\to\infty}\frac{\Ea^0[X_n\cdot e_1]}n\leq  \frac{b_1-b_3}{1+(1-p)[(b_1-b_3)+(r_3-r_1)]}.
$$
By the existence of deterministic limiting velocity and by bounded convergence, this is enough to yield the desired result.
\end{proof}

\section{Rare anomaly theorem}\label{sec:RareObs}

The goal of this section is to prove Theorem \ref{thm:MainThm} (of which Theorem \ref{thm:SpecialCase} is a special case) and Corollary \ref{cor:continuity}. To make the argument described in section \ref{subsec:ProofIdeas} precise, we first need to define what we mean by a ``cluster''. It turns out we only need to use a sort of directed cluster. Recall that $e_{\eta}$ and $e_{\theta}$ are distinct nearest neighbors of the origin with $e_{\eta}\neq- e_{\theta}$, and red sites satisfy a uniform ellipticity assumption that allows a walker to step from any red site $z$ to $z+e_{\theta}$ or $z+e_{\eta}$ with probability at least $\kappa$. Assume throughout this section (by relabling axes if necessary) that $\eta=2d-1$ and $\theta=2d$. Note that $e_{d}=-e_{2d}$. Let $\setM=\{e_{2d-1}\}\cup\{e_d\}\cup\{e_d+e_i:1\leq i \leq 2d,~i\notin\{d,2d\}\}$. If $d=3$, then $\setM$ is the ``cross'' above the origin, along with one vertex that is horizontally adjacent to the origin. Now in a given pre-environment $\pre$, define $C_x=C_x(\pre)$ to be the smallest subset of $\Z^d$ such that
\begin{itemize}
    \item If $x$ is red, then $x\in C_x$;
    \item If $y\in C_x$ and $z\in y+\setM$ is red, then $z\in C_x$. 
\end{itemize}
Equivalently, $C_x$ is the set of all sites $y$ that are red and can be reached by a \textit{red path} $\pi$ from $x$ to $y$ \textit{taking steps in $\setM$}; that is, a sequence $\pi=(x=\pi_0,\pi_1,\pi_2,\ldots,\pi_n=y)$ of all red vertices such that for all $1\leq k\leq n$, $(\pi_k-\pi_{k-1})\in\setM$ (note that such a path is not necessarily a nearest-neighbor path, as it can include certain diagonal steps). If $x$ is blue, then $C_x=\varnothing$.

\begin{figure}
    \centering
    \includegraphics[width=2in]{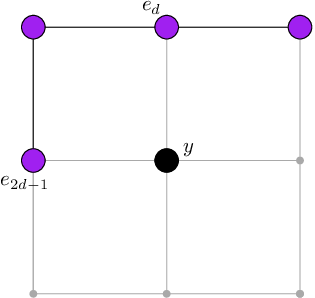}
    \hspace{1in}
    \includegraphics[width=3in]{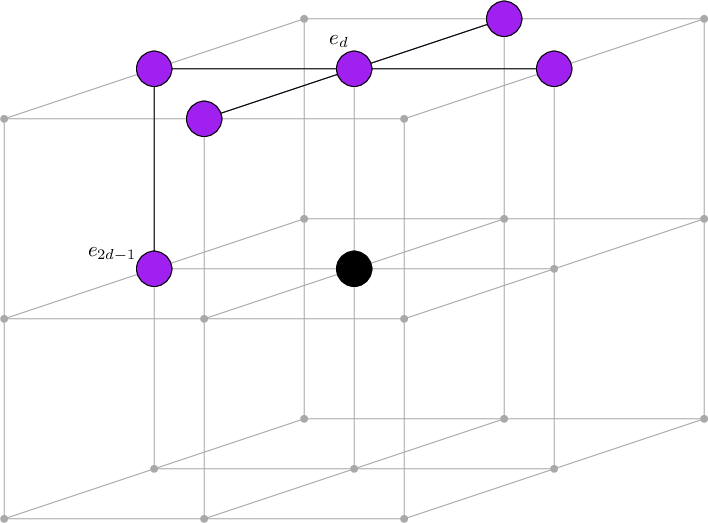}
    \caption{The set $\setM$ is depicted when $d=2$ on the left, and when $d=3$ on the right using an orthographic rendering. The origin, which is not in $\setM$, is black, and vertices in $\setM$ are purple.}
    \label{fig:SetM}
\end{figure}

The following proposition is the main ingredient in Theorem \ref{thm:MainThm}.

\begin{prop}\label{prop:MainProp}
    For all $\beta>0$, there exists $p^*\in(0,1)$ such that if $p>p^*$, then for all $x\in\Z^d$, for all $T\in\N$, and for every finite, nonempty set $C\subset\Z^d$, there exists $C'\subset C$ with $|C'|=|C|-1$ such that
    \begin{equation}\label{eqn:413}
        \Ea^0[N_x^T\mathbbm{1}_{\{C_x=C\}}]\leq\beta^{-1}\Ea^0[N_x^T\mathbbm{1}_{\{C_x=C'\}}].
    \end{equation}
\end{prop}

\begin{proof}
    Fix $x\in\Z^d$ and $\beta>0$. Let $C\subset\Z^d$ be finite and nonempty. We may freely assume that $\Pm(C_x=C)>0$, since otherwise \eqref{eqn:413} is trivial. Let $y$ be a vertex in $C$ maximizing $y\cdot e_d$, and of these (if there are more than one), let $y$ be a vertex maximizing $y\cdot e_{2d-1}$. We choose $y$ this way specifically to satisfy the following claim.
\begin{claim}\label{claim:422}
    If $C_x=C$, then every vertex in $y+\setM$ is blue. 
\end{claim}
To see that this is true, assume $C_x=C$ and note that every vertex in $z\in y+\setM$ has either $z\cdot e_d>y\cdot e_d$ or $z\cdot e_{2d-1}>y\cdot e_{2d-1}$, and therefore by the definition of $y$ cannot be in $C=C_x$. But since $y\in C_x$, the definition of $C_x$ says that any red vertex in $y+\setM$ must also be in $C_x$. Therefore, all vertices in $y+\setM$ are blue, proving Claim \ref{claim:422}.

Now let $C'=C\setminus\{y\}$. Rather than considering directly the event $\{C_x=C'\}$, we will consider a subset of this event. For a given pre-environment $\pre$, define the pre-environment $\pre^*$ by
\begin{equation}\notag
\pre^*(z)=\begin{cases}
    \pre(z)\quad\quad\quad&\text{if }z\neq y
    \\\notag
    \text{red}&\text{if }z=y
\end{cases}
\end{equation}
\begin{claim}\label{claim:435}
    If $C_x(\pre^*)=C$ and $\pre(y)=\text{blue}$, then $C_x(\pre)=C'$.
\end{claim}
    To see that $C_x(\pre)\subseteq C'$, note that switching $y$ from red to blue cannot create new red paths taking steps in $\setM$, and therefore $C_x(\pre)\subseteq C_x(\pre^*)=C$. Moreover, $y$ cannot be in $C_x(\pre)$ as it is blue, so $C_x(\pre)\subseteq C\setminus\{y\}=C'$.
    
    On the other hand, if $z\in C'$, then since $z\in C_x(\pre^*)$,  then there is a red path in $\pre^*$ from $x$ to $z$ taking steps in $\setM$; by Claim \ref{claim:422} and the fact that $z\neq y$, such a path cannot include $y$, so it is also a red path in $\pre$ from $x$ to $z$, implying $z\in C_x(\nu)$. Thus, $C'\subseteq C_x(\nu)$, proving Claim \ref{claim:435}. 
\begin{claim}\label{claim:426}
    $\Pm(C_x(\pre^*)=C,~\pre(y)=\text{blue})=\frac{p}{1-p}\Pm(C_x=C)$. 
\end{claim}
    To prove this claim, note that since $\pre^*$ is independent of $\pre(y)$, 
    \begin{equation}\label{eqn:443}
        \Pm(C_x(\pre)=C)=\Pm(C_x(\pre^*)=C)\Pm(y\text{ is red})=(1-p)\Pm(C_x(\pre^*)=C).
    \end{equation}
    Similarly,
    \begin{equation}\label{eqn:445}
        \Pm(C_x(\pre^*)=C,~\pre(y)=\text{blue})=\Pm(C_x(\pre^*)=C)\Pm(y\text{ is blue})=p\Pm(C_x(\pre^*)=C)
    \end{equation}
    Combining \eqref{eqn:443} and \eqref{eqn:445} proves Claim \ref{claim:426}.

For our next claim, we introduce a definition. For a set $C\subset \Z^d$, let $\partial C=(C+\setM)\setminus C$. 

\begin{claim}\label{claim:437}
    For a finite set $C\subset\Z^d$ with $\Pm(C_x=C)>0$, the event that $\{C_x=C\}$ is precisely the event that every vertex in $C$ is red and every vertex in $\partial C$ is blue. 
\end{claim}
 To prove this claim, note that $\Pm(C_x=C)>0$ if and only if $x\in C$ and every $z\in C$ can be reached from $x$ by a path in $C$ taking only steps in $\setM$. If every vertex in $C$ is red, this means every $z\in C$ can be reached from $x$ by a red path taking only steps in $\setM$, so $C\subseteq C_x$. But if every vertex in $\partial C$ is blue, then $C$ satisfies the conditions in the definition of $C_x$, and is therefore the smallest such set, meaning $C=C_x$. On the other hand, if $C_x=C$, it follows that every element of $C$ is red, and since any red vertex in $\partial C$ must also be in $C_x$, all vertices in $\partial C$ must be blue. This proves Claim \ref{claim:437}.

Now before making another claim, we will spend some time defining a coupling between three tagged environments and walks in them. 
Note that the event $\{C_x(\pre^*)=C\}$ is the event that all vertices in $C'$ are red in $\pre$ and all vertices in $\partial C$ are blue in $\pre$, with $y$ taking either color. Our goal is to have, in the first tagged environment, $C_x(\pre^*)\neq C$, in the second, $C_x(\pre^*)=C$ with $y$ blue, so that $C_x(\pre)=C'$, and in the third, $C_x(\pre)=C_x(\pre^*)=C$. We will define walks in all three environments, carefully coupling the last two walks with each other. The first only exists to allow us to create a tagged environment and walk drawn according to the annealed measure $\Pa^0$, and its dependence on the other two walks will not matter.

Consider, then, the measurable space given by the set
\begin{equation}\notag
\BigSpace=
\PRE\times\Omega\times\Sigma_{\mathcal{N}}^{C'\cup\partial C}\times\Sigma_{\mathcal{N}}^{2}\times[0,1]^{\N}\times[0,1]^{\N}\times\{1,2,3\}
\end{equation}
endowed with Borel (product) sigma algebra, and the measure $\widehat{P}$ that is the product of $\Pm(\cdot\given C_x(\pre^*)\neq C)$ on $\PRE\times\Omega$, $\mu_b^{\partial C}\times\mu_r^{C'}$ on $\Sigma_{\mathcal{N}}^{C'\cup\partial C}$, $\mu_b\times\mu_r$ on $\Sigma_{\mathcal{N}}^2$, the Borel-Lebesgue product measure on $[0,1]^{\N}\times[0,1]^{\N}$, and on $\{1,2,3\}$ the discrete probability measure that selects $1$ with probability $\Pm(C_x(\pre^*)\neq C)$,  selects $2$ with probability $\Pm(C_x(\pre^*)=C,~\pre(y)=\text{blue})$, and selects $3$ with probability $\Pm(C_x=C)$.
Say a canonical element of this measure space is the 8-tuple \begin{equation}\notag\tuple=(\pre_1,\omega_1,\omega^{C'\cup\partial C},\omega_b^y,\omega_r^y,{\bf U},{\bf V},A),\end{equation} where $A\in\{1,2,3\}$.

Let $\bomega_1:=(\pre_1,\omega_1)$. Define $\bomega_2=(\pre_2,\omega_2)$ to agree with $(\pre,\omega)$ everywhere except on $C\cup \partial C$, where $\pre_2$ is red on $C$ and blue on $\partial C$, and where $\omega_2^{C'\cup\partial C}=\omega^{C'\cup\partial C}$, and $\omega_2^y=\omega_b^y$. Define $(\pre_3,\omega_3)$ to agree with $(\pre_2,\omega_2)$ everywhere except at site $y$, where $\pre_3(y)=\text{red}$ and $\omega_3^y=\omega_r^y$.
The probability space is constructed so that by Claim \ref{claim:437} and the i.i.d. property of environments, the law of $(\pre_2,\omega_2)$ is the law under $\Pm$ of a tagged environment $\bomega$ conditioned to have $C_x(\pre^*)=C$ and $y$ blue, the law of $(\pre_3,\omega_3)$ is the law of a tagged environment conditioned to have $C_x=C$, and the law of $(\pre_A,\omega_A)$ is $\Pm$. 

The following lemma is standard. We state it here for completeness and to make sure our coupling construction that depends on it is clear.

\begin{lem}\label{lem:uniform}
    For all $a\in\Z^d$, there exists a map  $\varphi^a:\Omega\times[0,1]^{\N}\to(\Z^d)^{\N}$ taking a fixed environment $\omega$ and an infinite sequence ${\bf U}=(U_n)_{n=1}^{\infty}$ and outputting a walk $\varphi^a(\omega,{\bf U})=(X_n)_{n=0}^{\infty}$, such that
    \begin{enumerate}[(a)]
        \item For any fixed $\omega\in\Omega$, if ${\bf U}$ is a an i.i.d. sequence of uniform $(0,1)$ random variables, then the law of $\varphi^a(\omega,{\bf U})$ is $P_{\omega}^a$.
        \item If $\tau$ is a stopping time for $(X_n)_{n=0}^{\infty}$, then on the event $\tau<\infty$, the regular conditional law of $(U_n)_{n=\tau+1}^{\infty}$, conditioned on $(X_n)_{n=0}^{\tau}$, is almost surely the i.i.d. uniform law. 
    \end{enumerate}
\end{lem}

\begin{proof}[Proof of Lemma \ref{lem:uniform}]
    Let $X_0=a$. For $k\geq1$, let 
    \begin{equation}\notag i_k=\min\left\{1\leq i\leq 2d:\sum_{j=1}^i\omega^{X_{k-1}}(e_j)\leq U_k\right\},\end{equation}
    and let $\varphi^a(\omega,{\bf U})=(X_n)_{n=0}^{\infty}$, where $X_{n}=a+\sum_{k=1}^ne_{i_k}$. It is straightforward to check from this definition that (a) and (b) hold.
\end{proof}

We will define three coupled  walks ${\bf X}^1=(X_n^1)_{n=0}^{\infty}$, ${\bf X}^2=(X_n^2)_{n=0}^{\infty}$ and ${\bf X}^3=(X_n^3)_{n=0}^{\infty}$, drawn respectively according to $P_{\omega_1}^0$, $P_{\omega_2}^0$, and $P_{\omega_3}^0$. The only coupling that matters is between ${\bf X}^2$ and ${\bf X}^3$. The idea is that ${\bf X}^3$ steps with ${\bf X}^2$ until they hit site $y$ (the only site where transition probabilities in $\omega_2$ and $\omega_3$ can differ). After stepping away from $y$, the walk ${\bf X}^2$ ``freezes'', while ${\bf X}^3$ continues to move independently until hitting the neighbor of $y$ where ${\bf X}^2$ is waiting. Once the two walks meet, they recouple and move together. Each time they hit $y$, they follow the same protocol. If it is ever the case that ${\bf X}^3$ never hits the neighbor of $y$ that ${\bf X}^2$ stepped to, then the two walks run forever independently (we do not think of ${\bf X}^2$ as being frozen in this case).

Let us formalize the above description.
Let ${\bf X}^1$ be a walk started from the origin in environment $\omega_1$, using ${\bf U}$ to get its steps as in Lemma \ref{lem:uniform}. That is, ${\bf X}^1=\varphi^0(\omega_1,{\bf U})$. Similarly, let ${\bf X}^2=\varphi^0(\omega_2,{\bf U})$. By Lemma \ref{lem:uniform}, the law of ${\bf X}^1$, conditioned on $\omega_1$, is $P_{\omega_1}^0$, and similarly with ${\bf X}^2$.
We will define ${\bf X}^3$ in terms of ${\bf X}^2$, $\omega_3$, and ${\bf V}$. Let $\tau_1:=\inf\{n\geq0:X_n^2=y\}$ be the first hitting time of site $y$ by the walk ${\bf X}^2$, and recursively define $\tau_{k+1}:=\inf\{n>\tau_{k}:X_n^2=y\}$ as the time of the first return to $y$ after $\tau_k$. For $k$ such that $\tau_k<\infty$, define $y_k:=X_{\tau_k+1}^2$ to be the site ${\bf X}^2$ steps to after its $k$th visit to $y$.

Now if $\tau_1<\infty$, use the sequence ${\bf V}$ as in Lemma \ref{lem:uniform} to run a walk ${\bf Y}^1$ according to $P_{\omega_3}^x$ until time $b_1:=H_{y_1}$; that is, until the walk hits $y_1$, the location where ${\bf X}^2$ first stepped after hitting $y$. (Note that $b_1$ may be finite or infinite.) Now if $b_1,\ldots,b_k,\tau_{k+1}<\infty$, let $s_k:=\sum_{j=1}^kb_k$, and apply the function $\varphi^y$ from Lemma \ref{lem:uniform} to the environment $\omega_3$ and $(V_i)_{i=s_k+1}^{\infty}$ to run a walk ${\bf Y}^k$  until time $b_{k+1}:=H_{y_{k+1}}$, when the walk hits $y_{k+1}$.

Now define ${\bf X}^3:=(X_n^3)_{n=0}^{\infty}$ as the concatenation
\begin{equation}\label{eqn:317}
    (X_n^2)_{n=0}^{\tau_1-1}~\frown~{\bf Y}^1~\frown~(X_n^2)_{n=\tau_1+1}^{\tau_2-1}~\frown~{\bf Y}^2~\frown~(X_n^2)_{n=\tau_2+1}^{\tau_3-1}~\frown~{\bf Y}^3~\frown~\cdots,
\end{equation}
where $\frown$ denotes concatenation, and \eqref{eqn:317} is interpreted so that only terms that are defined are included. Note that our construction is such that with $\widehat{P}$-probability 1, if any sequence is defined, all previous sequences in the concatenation are defined and are finite, and conversely if the first $k$ sequences are finite, then the $(k+1)$st sequence is defined; this ensures that \eqref{eqn:317} almost surely gives us a single, well-defined, infinite sequence.

We have finished defining our coupling, and we now analyze it. We first define two conditional probability laws. Let $\tomega=(\pre_1,\omega_1,\omega^{C'\cup\partial C},\omega_b^y,\omega_r^y)$ (that is, the entries of $\tuple$ that determine the three tagged environments, but not the walks). Let $P_{\tomega}$ be the regular conditional law defined by $P_{\tomega}(\cdot)=\widehat{P}(\cdot\given \tomega)$. Let $\sigma({\bf V}^c)$ be the sigma field generated by every entry of $\tuple$ except ${\bf V}$, and let $Q(\cdot)=\widehat{P}(\cdot\given \sigma({\bf V}^c))$ be the law conditioned on $\sigma({\bf V}^c)$.

Our next three claims show that the walks we have defined obey the desired laws.
\begin{claim}\label{claim:524}
     It is $\widehat{P}$--almost surely the case that on the event that $\tau_1,\ldots,\tau_k<\infty$:
     \begin{itemize}
         \item the law of ${\bf Y}^1$ under $Q$ is the law under $P_{\omega_3}^x$ of a walk stopped at $H_{y_1}$; and
         \item for $j<k$, under the measure $Q$ conditioned on ${\bf Y}^1,\ldots,{\bf Y}^j$ all being finite;
         \begin{itemize}
             \item the law of ${\bf Y}^{j+1}$ is the law under $P_{\omega_3}^x$ of a walk stopped at $H_{y_{j+1}}$, and 
             \item ${\bf Y}^1,\ldots,{\bf Y}^{j+1}$ are independent. 
         \end{itemize}
     \end{itemize}  
\end{claim}
Claim \ref{claim:524} follows from induction by repeated application of Lemma \ref{lem:uniform}. 

\begin{claim}\label{claim:492}
    Conditioned on $\tomega$, each ${\bf X}^i$ is distributed according to $P_{\omega_i}^0$.
\end{claim}
For $i=1,2$, this is immediate from Lemma \ref{lem:uniform}.
 For $i=3$, it follows by induction with repeated application of the strong Markov property. Indeed, conditioned on the entries of $\tomega$, we know ${\bf X}^2$ follows $P_{\omega_3}^0$ until time $\tau_1$ , when the walk is at $y$ (since $\omega_2$ and $\omega_3$ agree except at $y$). Then ${\bf Y}^1$ follows the law of a walk in $\omega_3$ started at $y$, until it hits $y_1$. Then ${\bf X}^2$, started at time $\tau_1+1$, follows the law of a walk in $\omega_3$ started at $y_1$ until it hits $y$, and so on. This justifies Claim \ref{claim:492}.
\begin{claim}\label{claim:497}
    The law of $(\bomega_A,{\bf X}^A)$ under $\widehat{P}$ is $\Pa^0$. 
\end{claim}
To prove the claim, let $E\subset\bOmega$ and $F\subset(\Z^d)^{\N}$. Let $E_1=E\cap\{C_x(\pre^*)\neq C\}$, $E_2=E\cap\{C_x(\pre^*)=C,~\pre(y)=\text{blue}\}$, and $E_3=E\cap\{C_x(\pre)=C\}$. We have

\begin{align}\notag
    \widehat{P}(\bomega_A\in E_1,{\bf X}^A\in F)&=\widehat{P}(A=1,\bomega_1\in E_1,{\bf X}^1\in F)
    \\\notag
    &=\Pm(C_x(\pre^*)\neq C)\widehat{P}(\omega_1\in E_1,{\bf X}^1\in F)
    \\\notag
    &=\Pm(C_x(\pre^*)\neq C)\int_{\{\bomega_1\in E_1\}}\widehat{P}({\bf X}^1\in F~|~ \bomega_1)d\widehat{P}(\bomega_1\in\cdot)
    \\\notag
    &=\Pm(C_x(\pre^*)\neq C)\int_{\{\bomega_1\in E_1\}}P_{\omega_1}^0(F)d\widehat{P}(\bomega_1\in\cdot)
    \\\notag
    &=\Pm(C_x(\pre^*)\neq C)\int_{E_1}P_{\omega}^0(F)d\Pm(\cdot~|~ C_x(\pre^*)\neq C)
    \\\notag
    &=\int_{E_1}P_{\omega}^0(F)d\Pm
    \\\notag
    &=\Pa^0(E_1\times F).
\end{align}
Similarly, $\widehat{P}(\bomega_A\in E_2,{\bf X}^A\in F)=\Pa^0(E_2\times F)$ and $\widehat{P}(\bomega_A\in E_3,{\bf X}^A\in F)=\Pa^0(E_3\times F)$. Adding these together, we get $\widehat{P}(\bomega_A\in E,{\bf X}^A\in F)=\Pa^0(E\times F)$, which proves Claim \ref{claim:497}. 

We now begin a comparison of the walks ${\bf X}^2$ and ${\bf X}^3$, and in particular the numbers of times these walks hit $x$ before time $T$. In what follows, we let $a\wedge b:=\min(a,b)$. First we bound an expectation under the conditional measure $Q$.
\begin{claim}\label{claim:601}
With $\widehat{P}$-probability 1,
    \begin{align}\label{eqn:601}
    E_{Q}[N_x^T({\bf X}^3)] 
    &\leq N_x^T({\bf X}^2)+\sum_{k=1}^{T}\sum_{j=1}^{T}\mathbbm{1}_{\{\tau_k=j\}}\max_{i=1,\ldots,2d}E_{\omega_3}^y\left[N_x^{(T-j)\wedge H_{y+e_i}}\right]
\end{align}
\end{claim}
Examining \eqref{eqn:317}, note that if $\tau_k>T$ for some $k$, then all steps in ${\bf X}^3$ corresponding to ${\bf Y}^k$ occur after time $T$. We thus have, with $\widehat{P}$-probability 1, 
\begin{equation}\label{eqn:584}
    N_x^T({\bf X}^3)\leq N_x^T({\bf X}^2)+\sum_{k=1}^{T}\sum_{j=1}^{T}\sum_{i=1}^{2d}\mathbbm{1}_{\{\tau_k=j\}}\mathbbm{1}_{\{
    y_k
    =y+e_i\}}N_x^{T-j}({\bf Y}^k).
\end{equation}

Taking expectations under $Q$ in \eqref{eqn:584}, we almost surely get from Claim \ref{claim:524}:
\begin{align}\notag
    E_{Q}[N_x^T({\bf X}^3)]&\leq N_x^T({\bf X}^2)+\sum_{k=1}^{T}\sum_{j=1}^{T}\sum_{i=1}^{2d}\mathbbm{1}_{\{\tau_k=j\}}\mathbbm{1}_{\{
    y_k=y+e_i\}}E_Q\left[N_x^{T-j}({\bf Y}^k)\right].
    \\\notag
    &=N_x^T({\bf X}^2)+\sum_{k=1}^{T}\sum_{j=1}^{T}\sum_{i=1}^{2d}\mathbbm{1}_{\{\tau_k=j\}}\mathbbm{1}_{\{
    y_k=y+e_i\}}E_{\omega_3}^y\left[N_x^{(T-j)\wedge H_{y+e_i}}\right]
    \\\notag
    &\leq N_x^T({\bf X}^2)+\sum_{k=1}^{T}\sum_{j=1}^{T}\mathbbm{1}_{\{\tau_k=j\}}\max_{i=1,\ldots,2d}E_{\omega_3}^y\left[N_x^{(T-j)\wedge H_{y+e_i}}\right]
\end{align}
The last line is the assertion of Claim \ref{claim:601}. Our next claim puts a bound on the maximum inside the double sum.
\begin{claim}\label{claim:909}
    For each $1\leq i\leq 2d$, we have
    \begin{equation}\label{eqn:607}
        E_{\omega_3}^y\left[N_x^{(T-j)\wedge H_{y+e_i}}\right]
        \leq
        \max_{h=1,\ldots,2d}\frac{E_{\omega_3}^{y+e_h}\left[N_x^{(T-j-1)\wedge H_y}+\mathbbm1_{y=x}\right]}{\kappa^5}
    \end{equation}
\end{claim}

    To prove this claim, let $\siteb=y+e_i$. For a walk ${\bf X}=(X_n)_{n=0}^{\infty}$, define 
    \begin{equation}\notag
    T_j'=T_j'({\bf X}):=\inf\{n\geq T-j:y\notin\{
    X_{n-(T-j)+1},\ldots,X_{n-1},X_n\}\}
    \end{equation} 
    be the first time a walk has taken $T-j$ steps without hitting 
    $y$. We imagine that a clock allows the walk to take $T-j$ steps, but the clock resets 
    every time the walk hits $y$. The walk stops when the clock stops or when the walk hits $\siteb$. 
The stopping time $T_j'$ is useful because it allows us to use the strong Markov property every time the walk returns to $y$, but it is always at least $T-j$, so that we have 
\begin{equation}\label{eqn:625}
    E_{\omega_3}^y\left[N_x^{(T-j)\wedge H_{\siteb}}\right]
    \leq E_{\omega_3}^y\left[N_x^{T_j'\wedge H_{\siteb}}\right]
\end{equation}

By conditioning and then using the strong Markov property and the definition of $T_j'$, we can write\footnote{
Note that the event $\{\tilde{H}_y=T_j'\wedge H_{\siteb}\}$ is empty by definition.}

\begin{align}\notag
   E_{\omega_3}^y\left[N_x^{T_j'\wedge H_{\siteb}}\right]
    &=
    P_{\omega_3}^y(\tilde{H}_y> T_j'\wedge H_{\siteb})E_{\omega_3}^y\left[N_x^{T_j'\wedge H_{\siteb}}\given \tilde{H}_y> T_j'\wedge H_{\siteb}\right]
    \\\notag
    &\qquad\qquad +
    P_{\omega_3}^y(\tilde{H}_y<T_j'\wedge H_{\siteb})E_{\omega_3}^y\left[N_x^{T_j'\wedge H_{\siteb}}\given \tilde{H}_y<T_j'\wedge H_{\siteb}\right]
    \\\notag
    &=
    P_{\omega_3}^y(\tilde{H}_y> T_j'\wedge H_{\siteb})E_{\omega_3}^y\left[N_x^{T_j'\wedge H_{\siteb}}\given \tilde{H}_y> T_j'\wedge H_{\siteb}\right]
    \\\notag
    &\qquad\qquad +
    P_{\omega_3}^y(\tilde{H}_y<T_j'\wedge H_{\siteb})\left(E_{\omega_3}^y\left[N_x^{\tilde{H}_y-1}\given \tilde{H}_y<T_j'\wedge H_{\siteb}\right]+E_{\omega_3}^y\left[N_x^{T_j'\wedge H_{\siteb}}\right]\right)
    \\\notag
    &=E_{\omega_3}^y\left[N_x^{(T-j)\wedge H_{\siteb}\wedge (\tilde{H}_y-1)}\right]+P_{\omega_3}^y(\tilde{H}_y<T_j'\wedge H_{\siteb})E_{\omega_3}^y\left[N_x^{T_j'\wedge H_{\siteb}}\right]
\end{align}
To justify the last equality, consider the time $(T-j)\wedge H_b\wedge (\tilde{H}_y-1)$. On the event $\{\tilde{H}_y> T_j'\wedge H_{\siteb}\}$, this time is the same as $T_j'\wedge H_b$, while on the complementary event $\{\tilde{H}_y<T_j'\wedge H_{\siteb}\}$, it is the same as $\tilde{H}_y-1$.
Solving for $E_{\omega_3}^y\left[N_x^{T'\wedge H_b}\right]$ and noting that the complement of the event 
$\{\tilde{H}_y<T'\wedge H_b\}$ can also be expressed as $\{\tilde{H}_y> (T-j)\wedge H_b\}$
gives us
\begin{equation}\label{eqn:657}
    E_{\omega_3}^y\left[N_x^{T'\wedge H_b}\right]=\frac{E_{\omega_3}^y\left[N_x^{(T-j)\wedge H_b\wedge (\tilde{H}_y-1)}\right]}{P_{\omega_3}^y(\tilde{H}_y> (T-j)\wedge H_b)}.
\end{equation}
Now $P_{\omega_3}(\tilde{H}_y>(T-j)\wedge H_b)\geq P_{\omega_3}(H_b<\tilde{H}_y)$. If $b=y+e_{2d}$ or $b=y+e_{2d-1}$, then $P_{\omega_3}(H_b<\tilde{H}_y)\geq P_{\omega_3}(X_1=b)\geq \kappa$. If $b=y+e_d$, then $b$ can be reached by the path 
\begin{equation}\notag
(y+e_{2d-1}~,~y+e_{2d-1}+e_d~,~y+e_d)
\end{equation}
Any other neighbor of $y+e_i$ can be reached by the path
\begin{equation}\notag
(y~,~y+e_{2d-1}~,~y+e_{2d-1}+e_d~,~y+e_d~,~y+e_d+e_i~,~y+e_i).
\end{equation}
All sites in the above paths, except for $y$ itself $y+e_i$, are in $y+\setM$. They are therefore blue, and so the paths can be taken with probability at least $\kappa^4$ and $\kappa^5$, respectively. Thus for all neighbors $b$ of $y$, $P_{\omega_3}^y(H_b<\tilde{H}_y)\geq\kappa^5$.

\begin{figure}
    \centering
    \includegraphics[width=2in]{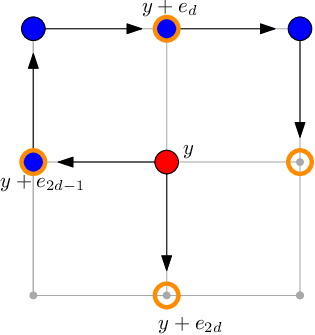}
    \hspace{1in}
    \includegraphics[width=3in]{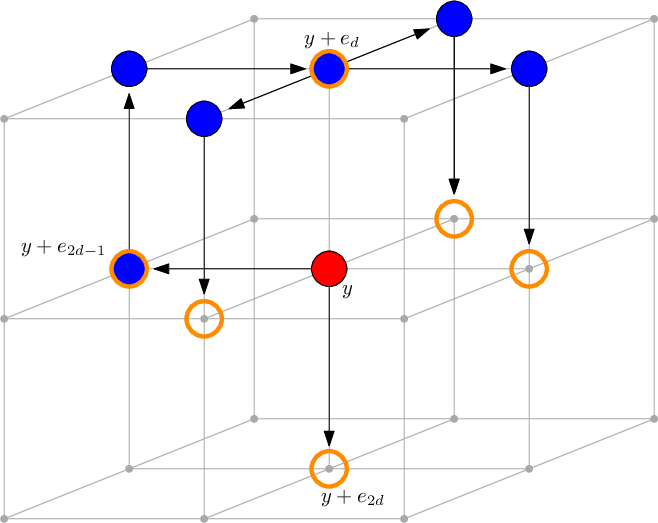}
    \caption{All nearest neighbors of $y$ (circled in orange) can be reached via paths of length no more than 5, where each step has probability at least $\kappa$ in $\omega_3$. The situations when $d=2$ and $d=3$ are depicted on the left and right, respectively. }
    \label{fig:Strategy}
\end{figure}

With these considerations, \eqref{eqn:657} gives us
\begin{equation}\notag
    E_{\omega_3}^y\left[N_x^{T'\wedge H_b}\right]\leq\frac{E_{\omega_3}^y\left[N_x^{(T-j)\wedge H_b\wedge (\tilde{H}_y-1)}\right]}{\kappa^5}.
\end{equation}
Now by considering where the walk first steps after visiting site $y$ (and throwing away $H_{\siteb}$), we can write
\begin{align}\notag
    E_{\omega_3}^y\left[N_x^{T'\wedge H_b}\right]&
    \leq
        \max_{h=1,\ldots,2d}\frac{E_{\omega_3}^{y+e_h}\left[N_x^{(T-j-1)\wedge (H_y-1)}+\mathbbm1_{y=x}\right]}{\kappa^5},
\end{align}
which is \eqref{eqn:607}. This proves Claim \ref{claim:909}.

Our next claim is an almost-sure lower bound for $E_{\tomega}[N_x^T({\bf X}^3)]$, the expectation under $P_{\tomega}$.
\begin{claim}\label{claim:1058} With probability 1,
\begin{equation}\label{eqn:729}
    E_{\tomega}[N_x^T({\bf X}^3)]\leq E_{\tomega}[N_x^T({\bf X}^2)]+ \sum_{k=1}^{T}\sum_{j=1}^{T}P_{\tomega}(\tau_k=j)\max_{h=1,\ldots,2d}\frac{E_{\omega_2}^{y+e_h}\left[N_x^{(T-j-1)\wedge (H_y-1)}+\mathbbm1_{y=x}\right]}{\kappa^5}.
\end{equation}
\end{claim}

Putting \eqref{eqn:607} from Claim \ref{claim:909} into \eqref{eqn:601} from Claim \ref{claim:601}, we get
\begin{align}\notag
    E_Q[N_x^T({\bf X}^3)]
&\leq 
N_x^T({\bf X}^2)+\sum_{k=1}^{T}\sum_{j=1}^{T}\mathbbm{1}_{\{\tau_k=j\}}
\max_{h=1,\ldots,2d}\frac{E_{\omega_3}^{y+e_h}\left[N_x^{(T-j-1)\wedge (H_y-1)}+\mathbbm1_{y=x}\right]}{\kappa^5}
\\\label{eqn:687}
&= 
N_x^T({\bf X}^2)+\sum_{k=1}^{T}\sum_{j=1}^{T}\mathbbm{1}_{\{\tau_k=j\}}
\max_{h=1,\ldots,2d}\frac{E_{\omega_2}^{y+e_h}\left[N_x^{(T-j-1)\wedge (H_y-1)}+\mathbbm1_{y=x}\right]}{\kappa^5}.
\end{align}
The equality comes from the fact that $\omega_2$ and $\omega_3$ agree away from site $y$.

Now because $\sigma(\tomega)\subset\sigma({\bf V}^c)$, we almost surely have $E_{\tomega}[E_Q[\cdot]]=E_{\tomega}[\cdot]$. Therefore, taking expectations in \eqref{eqn:687} proves Claim \ref{claim:1058}.

We now claim an almost-sure upper bound for $E_{\tomega}[N_x^T({\bf X}^2)]$.
\begin{claim}\label{claim:1075} With probability 1,
\begin{align}\label{eqn:723}
    E_{\tomega}[N_x^T({\bf X}^2)]\geq
    \sum_{k=1}^{T}\sum_{j=1}^{T}P_{\tomega}(\tau_k=j)
\max_{h=1,\ldots,2d}\kappa E_{\omega_2}^{y+e_h}\left[N_x^{(T-j-1)\wedge (H_y-1)}+\mathbbm1_{y=x}\right].
\end{align}
\end{claim}
We can sort $N_x^T({\bf X}^2)$ in terms of the visits of ${\bf X}^2$ to site $y$, throwing away any visits to $x$ before the first visit to $y$, to get
\begin{align}\label{eqn:714}
    N_x^T({\bf X}^2)\geq\sum_{k=1}^T\sum_{j=1}^T\mathbbm{1}_{\tau_k=j}\#\{j\leq n\leq T\wedge (\tau_{k+1}-1):X_n=x\}
\end{align}
Taking expectations in \eqref{eqn:714} with respect to $P_{\tomega}$, we get (by Claim \ref{claim:492} and the strong Markov property):
\begin{align}\label{eqn:718}
    E_{\tomega}[N_x^T({\bf X}^2)]\geq\sum_{k=1}^T\sum_{j=1}^TP_{\tomega}(\tau_k=j)
    E_{\omega_2}^y[N_x^{(T-j)\wedge(H_y-1)}]
\end{align}
Now because $y$ is blue in $\omega_2$, the a walk has at least $P_{\omega_2}^y$-probability $\kappa$ of stepping to each of its neighbors on the first step, including the neighbor that maximizes $E_{\omega_2}^{y+e_h}[N_x^{(T-j-1)\wedge(H_y-1)}]$. By restricting to the event that the walk from $y$ steps to this site and applying the strong Markov property, we almost surely get \eqref{eqn:723}, which proves Claim \ref{claim:1075}.

Our next (and final) claim combines the lower bound and the upper bound from the previous two claims, and translates the comparison into a statement about expectations under the original annealed measure. 
\begin{claim}\label{claim:1113}
\begin{equation}\label{eqn:1114}
    \Ea^0[N_x^T\mathbbm{1}_{\{C_x=C\}}]\leq\frac{1-p}p\left(1+\frac1{\kappa^6}\right)\Ea^0\left[N_x^T\mathbbm{1}_{\{C_x=C'\}}\right].
\end{equation}
\end{claim}
Combining \eqref{eqn:729} from Claim \ref{claim:1058} with \eqref{eqn:723} from Claim \ref{claim:1075} gives us, almost surely,
\begin{equation}\label{eqn:734}
    E_{\tomega}[N_x^T({\bf X}^3)]\leq\left(1+\frac{1}{\kappa^6}\right)E_{\tomega}[N_x^T({\bf X}^2)].
\end{equation}
Now
\begin{align}\notag
    \Ea^0\left[N_x^T\mathbbm{1}_{\{C_x=C\}}\right]
    &=
    \widehat{E}\left[N_x^T({\bf X}^A)\mathbbm{1}_{\{C_x({\bomega_A})=C\}}\right]
    \\\notag
    &=\widehat{E}\left[N_x^T({\bf X}^3)\mathbbm{1}_{\{A=3\}}\right]
    \\\notag
    &=\widehat{E}\left[N_x^T({\bf X}^3)\right]\Pm(C_x=C)
    \\\label{eqn:748}
    &=\widehat{E}\left[ E_{\tomega}[N_x^T({\bf X}^3)]\right]\Pm(C_x=C)
\end{align}
The first equality comes from Claim \ref{claim:497}. The third comes from independence of $A$ from ${\bf X}^3$.

Now by applying Claim \ref{claim:435}, arguing as above, and then applying Claim \ref{claim:426}, we get
\begin{align}\notag
    \Ea^0\left[N_x^T\mathbbm{1}_{\{C_x=C'\}}\right]
    &\geq
    \Ea^0\left[N_x^T\mathbbm{1}_{\{C_x(\pre^*)=C,~\pre(y)=\text{blue}\}}\right]
    \\\notag
    &=
    \widehat{E}\left[ E_{\tomega}[N_x^T({\bf X}^2)]\right]\Pm(C_x(\pre^*)=C,~\pre(y)=\text{blue})
    \\\label{eqn:759}
    &=
    \widehat{E}\left[ E_{\tomega}[N_x^T({\bf X}^2)]\right]\frac{p}{1-p}\Pm(C_x=C)
\end{align}
Now \eqref{eqn:734}, \eqref{eqn:748}, and \eqref{eqn:759} combine to yield \eqref{eqn:1114}, which proves Claim \ref{claim:1113}.

Now
for any $\beta>0$ and fixed $\kappa$, we may choose $p^*$ close enough to 1 that for any $p>p^*$, we have $\frac{1-p}{p}\left(1+\frac{1}{\kappa^6}\right)\leq\beta^{-1}$. Putting this into  Claim \ref{claim:1113} finishes the proof of Proposition \ref{prop:MainProp}.
\end{proof}

Our next goal is to translate Proposition \ref{prop:MainProp} into a comparison between the truncated Red's and Blue's functions. For this, we require a bound on the number of possible sets $C_x$ of various sizes, which we get via a standard lattice-animal-type argument. We base our proof on the argument for Equation 4.24 in \cite{Grimmett1999}), writing it out in full because $C_x$ may not be a ``lattice animal'' in the traditional sense, because writing the argument here allows us to get a more optimal bound for our situation, and because the argument is short and in fact simpler in our setting than in that of \cite{Grimmett1999}.

\begin{lem}\label{lem:latticeanimalbound}
    The number of sets $C$ such that $|C|=n$ and $\Pm(C_0=C)>0$ is no more than $(7d)^n$.
\end{lem}

\begin{proof}
    Let $\mathcal{C}_0^{n,b}$ be the set of sets $C\subset\Z^d$ such that $|C|=n$, $|\partial C|=b$, and $\Pm(C_0=C)>0$, and let $a_{n,b}=|\mathcal{C}_0^{n,b}|$. By Claim \ref{claim:437} from the proof of Proposition \ref{prop:MainProp}, we have, for fixed $n$,
    \begin{equation}\label{eqn:1138}
        \sum_{b}a_{n,b}p^n(1-p)^b=\Pm(|C|=n)\leq1.
    \end{equation}
    Now because $|\setM|=2d$, we must have $|\partial C_0|\leq2d|C_0|$, so $\mathcal{C}_0^{n,b}$ is empty whenever $b>2dn$. From \eqref{eqn:1138}, we therefore have
    \begin{equation}\notag
        \sum_{b}a_{n,b}p^n(1-p)^{2dn}\leq1,
    \end{equation}
    which gives us
    \begin{equation}\notag
        \sum_{b}a_{n,b}\leq\left(p(1-p)^{2d}\right)^{-n}.
    \end{equation}
    By optimizing\footnote{The maximum value of $p(1-p)^{2d}$ as $p$ ranges from 0 to 1 is increasing in $d$ and is asymptotically $\frac1{2ed}\approx\frac{1}{5.44d}$ as $d$ goes to infinity. For $d\geq3$ we may take 6 in place of 7 here. For $d=2$ we may take 6.11, but we use 7 in our proofs because it is an integer that works for all $d\geq2$. Nevertheless, a number somewhat lower than 6.5 is needed to get the ``13'' in Remark \ref{rem:exactbound}.} over all $p\in[0,1]$, we can always find a value of $p$ such that $p(1-p)^{2d}\geq \frac{1}{7d}$, so $\sum_{b}a_{n,b}\leq(7d)^n$.
\end{proof}

We can now use our bound on clusters to get a comparison between truncated Red's and Blue's functions.

\begin{prop}\label{prop2}
    For every $\varepsilon>0$, there exists $p^*<1$, depending only on $\varepsilon$, $d$, and $\kappa$, such that if $p>p^*$, then 
    \begin{equation}\label{eqn:979}
        \Ea^0[N_x^T\mathbbm{1}_{\{x\text{ is red}\}}]
        \leq
        \varepsilon\Ea^0[N_x^T\mathbbm{1}_{\{x\text{ is blue}\}}]
    \end{equation}
\end{prop}

\begin{proof}
    
Let $\varepsilon>0$. Choose $\beta=7d\left(1+\frac{1}{\varepsilon}\right)$. By Proposition \ref{prop:MainProp} and Lemma \ref{lem:condC}, we may choose $p^*<1$ large enough that for $p>p^*$, condition (C) is satisfied and for all $x\in\Z^d$, $T\in\N$, and for every finite, nonempty set $C\subset\Z^d$, there exists $C'\subset C$ with $|C'|=|C|-1$ such that 
    \begin{equation}\notag
        \Ea^0\left[N_x^T\mathbbm{1}_{\{C_x=C'\}}\right]
        \geq \beta
        \Ea^0\left[N_x^T\mathbbm{1}_{\{C_x=C\}}\right].
    \end{equation}
By induction, this implies that for every finite, nonempty set $C\subset\Z^d$, 
\begin{align}\notag
    \Ea^0\left[N_x^T\mathbbm{1}_{\{C_x=C\}}\right]
    &\leq
    \beta^{-|C|}\Ea^0\left[N_x^T\mathbbm{1}_{\{C_x=\varnothing\}}\right]
    \\\label{eqn:784}
    &=\beta^{-|C|}\Ea^0\left[N_x^T\mathbbm{1}_{\{x\text{ is blue}\}}\right].
\end{align}
 Let $x\in\Z^d$. By Lemma \ref{lem:latticeanimalbound} (and shift-invariance of $\Pm$), there are for each $k\geq1$ at most $(7d)^k$ sets $C$ such that $|C|=k$ and $\Pm(C_x=C)>0$. Let $\mathcal{C}_x^k$ be the collection of all such sets. Choose also $p^*$ large enough that for $p>p^*$, $|C_x|<\infty$ for all $x\in\Z^d$ (a standard union-bound argument shows that if $p>1-\frac{1}{2d}$, then the probability of there being a red path of length $k$ from $x$ taking steps in $\setM$ decays as $k\to\infty$, which implies $|C_x|$ must be finite). We then have
    \begin{align}\notag
        \Ea^0[N_x^T\mathbbm{1}_{\{x\text{ is red}\}}]
        &=
        \sum_{k=1}^{\infty}\sum_{C\in\mathcal{C}_x^k}\Ea^0[N_x^T\mathbbm{1}_{\{C_x=C\}}]
        \\\notag
        &\overset{\eqref{eqn:784}}{\leq}
        \sum_{k=1}^{\infty}(7d)^{k}\beta^{-k}\Ea^0[N_x^T\mathbbm{1}_{\{x\text{ is blue}\}}]
        \\\notag
        &=\Ea^0[N_x^T\mathbbm{1}_{\{x\text{ is blue}\}}]\sum_{k=1}^{\infty}\left(1+\frac{1}{\varepsilon}\right)^{-k}
        \\\notag
        &=\varepsilon\Ea^0[N_x^T\mathbbm{1}_{\{x\text{ is blue}\}}],
    \end{align}
which completes the proof.
\end{proof}

We are ready to prove our main theorem.

\begin{proof}[Proof of Theorem \ref{thm:MainThm}]
Fix $\varepsilon>0$. By Proposition \ref{prop2}, we may choose $p^*\geq p_c(d)$ satisfying \eqref{eqn:395c}. Then for $p>p^*$, by Lemma \ref{lem:condC} and Corollary \ref{cor:DirectionalTransience}, we have almost surely
    \begin{align}\notag
    v_1^b-\varepsilon\leq \lim_{n\to\infty}\frac{X_n\cdot\vec{u}}{n}
    &\leq v_2^b+\varepsilon,
\end{align}
which is what we wanted to prove.
\end{proof}

\begin{proof}[Proof of Theorem \ref{thm:SpecialCase}]
    This is a special case of Theorem \ref{thm:MainThm} where $\mu_b$ and $\mu_r$ are each a Dirac mass at one probability vector. Here $\kappa=.2$, $\theta=4$, and $\eta=3$ if $\delta<\frac14$; otherwise $\eta=1$.
\end{proof}

\begin{proof}[Proof of Corollary \ref{cor:continuity}]

Suppose $\mu$ is a Dirac mass at $\xi\in\Sigma_{\mathcal{N}}^{\kappa,\eta,\theta}$. From the strong law of large numbers, $\vel(\mu)$ is the Dirac mass at the vector $\left(\xi_i-\xi_{d+i}\right)_{i=1}^d$. Call this vector $\vect{v}(\xi)$. Indeed, for any $\xi'\in\Sigma_{\mathcal{N}}^{\kappa,\eta,\theta}$, let $\vect{v}(\xi'):=(\xi_i-\xi_{d+i})_{i=1}^d$ (thus, we have $\vel(\delta_{\xi'})=\delta_{\vect{v}(\xi')}$). 

Let $\mu_n$ be a sequence of measures converging weakly to $\mu$, and let $\varepsilon>0$. Let $\Pa_n^0$ be the annealed measure associated with the RWRE where all transition probability vectors are drawn according to $\mu_n$.
It suffices to show that for large enough $n$, the support of $\vel(\mu_n)$ is contained in an $\varepsilon$-neighborhood of the vector $\vect{v}(\xi)$, or equivalently that for large enough $n$,
\begin{equation}
    \label{eqn:1208}
    \left|\lim_{n\to\infty}\frac{X_n}{n}-\vect{v}(\xi)\right|<\varepsilon,~\Pa_n^0-a.s.
\end{equation}

     Choose $\delta :=\frac{\varepsilon}{2\sqrt{d}}$. Then for any vector $\vec{w}\in\R^d$, if $|\vec{w}\cdot e_i-\vect{v}(\xi)\cdot e_i|<\delta$ for all $i$, then $|\vec{w}-\vect{v}(\xi)|<{\varepsilon}$. Now choose $\varepsilon'$ such that if $|\xi'-\xi|<\varepsilon'$, then $|\vect{v}(\xi')\cdot e_i-\vect{v}(\xi)\cdot e_i|<\frac{\delta}{2}$ for all $i$ (this can be done by the continuity of $\vect{v}$ as a function of $\xi'$). Now for $r>0$ let $B_{r}(\xi):=\{\xi'\in\Sigma_{\mathcal{N}}^{\kappa,\eta,\theta}:|\xi'-\xi|<\varepsilon'\}$. 
    Assume $\varepsilon'$ is also chosen small enough that for $\xi'\in B_{\varepsilon'}(\xi)$, every component of $\xi'$ is at least $\kappa'$ for some $\kappa'>0$---we can do this because every component of $\xi$ is nonzero.
    
    By the fact that $\mu_n$ converges weakly to $\mu$ and the fact that $\mu(B_{\varepsilon'}(\xi))=1$, we have $\lim_{n\to\infty}\mu_n(B_{\varepsilon'}(\xi))=1$. 
    Now for each $n$ with $\mu_n(B_{\varepsilon'}(\xi))\in(0,1)$, define the measures $\mu_{n,b}$ and $\mu_{n,r}$ on $\Sigma_{\mathcal{N}}^{\kappa,\eta,\theta}$ by $\mu_{n,b}(A)=\frac{\mu_{n}(A\cap B_{\varepsilon'}(\xi))}{\mu_n(B_{\varepsilon'}(\xi))}$ and $\mu_{n,r}(A)=\frac{\mu_{n}(A\cap B_{\varepsilon'}(\xi)^c)}{\mu_n(B_{\varepsilon'}(\xi)^c)}$, and let $\Pa_{n,b}^0$ and $\Pa_{n,r}^0$ be the associated annealed measures.  

    Each vector $\xi'$ in the support of $\mu_{n,b}$ has every component at least $\kappa'$ and has $|\vect{v}(\xi')\cdot  e_i-\vect{v}(\xi)\cdot e_i|<\frac{\delta}{2}$ for all $i$. Moreover, every vector in the support of $\mu_{r,n}$ has $e_{\theta}$ and $e_{\eta}$ components at least $\kappa$. Therefore, for each $1\leq i\leq d$, $\mu_n$ satisfies the assumptions of Theorem \ref{thm:MainThm} with $\vec{u}=e_i$, $v_1=\vect{v}(\xi)\cdot e_i-\frac{\delta}{2}$, and $v_2=\vect{v}(\xi)\cdot e_i+\frac{\delta}{2}$. Since $\mu_n(B_{\varepsilon'}(\xi))$ approaches 1, we may conclude from Theorem \ref{thm:MainThm} 
    that for large enough $n$,
    $$
    \vect{v}(\xi)\cdot e_i-\delta<\lim_{n\to\infty}\frac{X_n\cdot e_i}{n}<\vect{v}(\xi)\cdot e_i+\delta,~\Pa_{n}^0-a.s.
    $$
    for all $1\leq i\leq d$, whence we may conclude \eqref{eqn:1208}.
\end{proof}

\section{A strongly mixing counterexample}\label{sec:Counterexample}

In this section, we show the significance of the i.i.d. assumption in Theorem \ref{thm:MainThm} by providing a counterexample to the statement---in fact, a counterexample to the statement of Theorem \ref{thm:SpecialCase}---with the i.i.d. assumption removed. For $\beta>0$ and $p\in(0,1)$, we will present a $\delta>0$ and a measure $\Pm$ on environments where red and blue sites are as in Theorem \ref{thm:SpecialCase}, where every site is blue with probability at least $p$, and where $\lim_{n\to\infty}\frac{X_n\cdot e_1}{n}<\beta$.

\begin{exmp}\label{exmp:Counterexample}
    Let $(U_x)_{x\in\Z^d}$ be an i.i.d. collection of uniform $(0,1)$ random variables. Let $\varepsilon>0$. For each $x=(a,b)\in\Z^d$, let $x^*$ be the closest vertex directly above $x$ whose associated uniform random variable is less than $\varepsilon$. In other words,  $x^*=(a,b^*)$, where $b^*=\inf\{b'\geq b:U_{(a,b')}<\varepsilon\}$. For each $x$, we have $U_{x^*}<\varepsilon$ by definition. Let $x$ be red if in fact $U_{x^*}<\varepsilon^2$; otherwise, let $x$ be blue, where red and blue sites are as in Theorem \ref{thm:SpecialCase}. Let $\Pm$ be the pushforward measure on environments.
\end{exmp}

\begin{prop}
    Let $p\in(0,1)$ and $\beta>0$. For appropriately chosen $\delta>0$ and $\varepsilon>0$, the measure $\Pm$ described in Example \ref{exmp:Counterexample} is such that sites are blue with probability at least $p$, and yet $\lim_{n\to\infty}\frac{X_n\cdot e_1}{n}<\beta$.
\end{prop}

\begin{proof}
    We will defer our choice of $\varepsilon$ to the end of the proof. For now, assume $\varepsilon<\frac12$, and let $\delta=\varepsilon^2$.
    Note that for each $x\in\Z^d$, $x^*$ is independent of $U_{x^*}$, which is uniformly distributed on $(0,\varepsilon)$. Now let $x_*$ be the nearest site directly \textit{below} $x$ such that $U_{x_*}<\varepsilon$. The distance from $x$ to $x_*$ and one plus the distance from $x$ to $x^*$ are independent geometric random variables with parameter $\varepsilon$. Hence there exists $C_1$ not depending on $\varepsilon$ (except through the fact that $\varepsilon<\frac12$) such that with probability at least $C_1$, $|x^*-x|$ and $|x-x_*|$ are both greater than $\frac1{\varepsilon}$. With probability at least $C_1\varepsilon$, then, $|x^*-x|,~|x-x_*|<\frac{1}{\varepsilon}$ and $U_{x^*}<\varepsilon^2$. On this event, all sites strictly between $x_*$ and $x^*$, as well as $x^*$, are red. 

    Now under the measure $\Pa^0$, for each $k\geq1$ let $T_k=\inf\{n\in\N:X_n\cdot e_1\geq k\}$ be the $k$th time the walk enters a fresh half plane to the right of the origin. 
    Let $\mathcal{F}_k$ be the sigma field generated by $(\omega^x)_{x\cdot e_1<k}$ as well as $(X_n)_{n=0}^{T_k}$. We consider $\Ea^0[T_{k+1}-T_k\given \mathcal{F}_k]$. 

    Say $y=X_{T_k}$. Let $E_k$ be the event that $|y^*-y|,~|y-y_*|>\frac{1}{\varepsilon}$ and $U_{y^*}<\varepsilon^2$. By shift invariance and the fact that $E_k$ depends only on the environment to the right of sites that have been hit before time $T_k$, the event $E_k$ is independent of $\mathcal{F}_k$ and $P(E_k)>C_1\delta$. On $E_k$, hitting the half-plane $\{x\cdot e_1\geq k+1\}$ requires the walk to either pass through a red site, which typically takes at least $O(\frac{1}{\varepsilon^2})$ steps to happen, or to move around a wall of $\frac{1}{\varepsilon}$ vertices in each direction. All vertical steps are up with probability $\frac14$ and down with probability $\frac14$, regardless of site color. Therefore, reaching a site with vertical position that differs from that of $y$ by $\frac{1}{\varepsilon}$ or more typically takes $O(\frac{1}{\varepsilon^2})$ steps. In fact, on the event $E_k$, with $P(\cdot\given X_1,\ldots,X_{T_k})$-probability at least $C_2$, it takes at more than $\frac{1}{\varepsilon^2}$ steps to hit the half-plane $\{x\cdot e_1\geq k+1\}$. Therefore, with probability 1,
    \begin{align}\notag
    \Ea^0[T_{k+1}-T_k\given \mathcal{F}_k]&\geq 
    C_1\varepsilon\Ea^0[T_{k+1}-T_k\given \mathcal{F}_k]
    \\\notag
    &\geq
    C_1\varepsilon C_2\frac{1}{\varepsilon^2}
    \\\notag
    &=\frac{C_1C_2}{\varepsilon}.
    \end{align}
    It follows that
    \begin{equation}\notag
    \Ea^0[T_{k+1}-T_k]\geq C_1\varepsilon C_2\frac{1}{\varepsilon^2}=\frac{C_1C_2}{\varepsilon},
    \end{equation}
    whence it is standard (see, e.g., \cite[Lemma 2.1.17]{Zeitouni2004}) to get
    \begin{equation}\notag
    \lim_{n\to\infty}\frac{X_n\cdot e_1}{n}\leq \frac{\varepsilon}{C_1C_2}.
    \end{equation}
    Now if we choose $\varepsilon<\min(1-p,C_1C_2\beta)$, then sites are blue with probability at least $p$, and $\lim_{n\to\infty}\frac{X_n\cdot e_1}{n}<\beta$.
\end{proof}

We now show that the measure $\Pm$ described in Example \ref{exmp:Counterexample} has good mixing properties.

\begin{prop}
    Under the measure $\Pm$ described in Example \ref{exmp:Counterexample}, the environment is $\psi$-mixing along vertical strips, and the vertical strips are i.i.d.  
\end{prop}

\begin{proof}
    That the vertical strips are i.i.d. follows immediately from the construction of $\Pm$. We now show that the environment is $\psi$-mixing along a vertical strip $\{(a,y):y\in\Z\}$. That the environment is strictly stationary along this strip follows, again, immediately from the construction. Now without loss of generality let $a=0$. Let $\mathcal{A}$ be the sigma field generated by $(\omega^{(a,y)})_{y\leq0}$, and for $n\geq1$, let $\mathcal{B}_n$ be the sigma field generated by $(\omega^{(a,y)})_{y\geq n}$. Recall that $\psi$-mixing means
    \begin{equation}\label{eqn:1042}
        \lim_{n\to\infty}\sup_{\substack{A\in\mathcal{A},~B\in\mathcal{B}_n \\\notag \Pm(A),\Pm(B)>0}}\left[\frac{\Pm(A\cap B)}{\Pm(A)\Pm(B)}-1\right]=0
    \end{equation}
    Let $P$ be the measure on the uniform variables used to determine $\omega$. By slight abuse of notation, if $E$ is an event in the sigma field on environments, we use $P(E)$ to denote the probability that the configuration of uniform random variables is such that $\omega\in E$. Using this convention, \eqref{eqn:1042} is equivalent to the same equation with $\Pm$ replaced by $P$ throughout.

    Now let $A\in\mathcal{A}$, and $B\in\mathcal{B}$, both with positive probability. Let $\mathcal{C}$ and $\mathcal{D}$ denote, respectively, the sigma fields generated by the uniform variables assigned to sites strictly below $(0,n)$ and above $(0,n)$, inclusive.
    
    By construction, $B\in\mathcal{D}$.
    Let $M$ be the event that $(0,0)^*<n$. Equivalently, $M$ is the event that $U_{(0,b)}<\varepsilon$ for some $b<n$. On the event $M$, the environment at and below $(0,0)$ is completely determined by the values of uniform random variables assigned to sites strictly below $(0,n)$. Therefore, $A\cap M\in\mathcal{C}$, so that $B$ is independent of $A\cap M$. Moreover, $A$ is independent of $M$, since the value of $U_{(0,0)^*}$ is independent of $(0,0)^*$. We therefore have
    \begin{align}\notag
        P(A\cap M\cap B)
        &=P(A\cap M)P(B)
        \\\label{eqn:1054}
        &=P(A)P(M)P(B)
    \end{align}
    so that
    \begin{align}\notag
        P(A\cap B)&=P(A\cap B\cap M^c)+P(A\cap M\cap B)\\\label{eqn:1058}
        &\geq P(A)P(B)P(M).
    \end{align}
This gives us a lower bound on $P(A\cap B)$ in terms of $P(A)P(B)$. We now seek an upper bound. We have
\begin{align}\label{eqn:1061}
    P(A\cap B\cap M^c)&=P(A\cap B\cap M^c\cap\{U_{(0,n)^*}<\varepsilon^2\})+P(A\cap B\cap M^c\cap\{U_{(0,n)^*}>\varepsilon^2\})
\end{align}
We examine the first term on the right side of \eqref{eqn:1061}. Note that $M$ is independent of $\mathcal{D}$, and in particular $M$ is independent of $B$ and of the event that $U_{(0,n)^*}<\varepsilon$. Note further that on the event $M^c$, $(0,0)^*=(0,n)^*$, but $A$ depends on $\mathcal{D}$ only through $U_{(0,n)^*}$, and in particular whether it is less than or greater than $\varepsilon^2$. We therefore have
\begin{align}\notag
    P(A\cap B\cap M^c\cap\{U_{(0,n)^*}<\varepsilon^2\})
    &=P(U_{(0,n)^*}<\varepsilon^2)P(B\given U_{(0,n)^*}<\varepsilon^2)
    \\\notag
    &\qquad\qquad \times P(A\cap M^c\given B\cap\{U_{(0,n)^*}<\varepsilon^2\})
    \\\notag
    &=\varepsilon P(B\given U_{(0,n)^*}<\varepsilon^2)P(M^c)P(A\given  M^c\cap B\cap\{U_{(0,n)^*}<\varepsilon^2\})
    \\\notag
    &=\varepsilon P(B\given U_{(0,n)^*}<\varepsilon^2)P(M^c)P(A\given U_{(0,0)^*}<\varepsilon^2).
\end{align}
On the other hand, 
\begin{align}\notag
    P(A\cap B\cap M\cap\{U_{(0,n)^*}<\varepsilon^2\})&\geq
    P(A\cap B\cap M\cap\{U_{(0,n)^*}<\varepsilon^2\}\cap\{U_{(0,0)^*}<\varepsilon^2\})
    \\\notag
    &=P(M)P(U_{(0,0)^*}<\varepsilon^2)P(U_{(0,n)^*}<\varepsilon^2)
    \\\notag
    &\qquad\qquad \times P(A\given U_{(0,0)^*}<\varepsilon^2)P(B\given U_{(0,n)^*}<\varepsilon^2)
    \\\notag
    &=P(M)\varepsilon^2P(A\given U_{(0,0)^*}<\varepsilon^2)P(B\given U_{(0,n)^*}<\varepsilon^2)
    \\\notag
    &= \frac{\varepsilon P(M)}{P(M^c)}P(A\cap B\cap M^c\cap\{U_{(0,n)^*}<\varepsilon^2\})
\end{align}
An analogous calculation can be done for the second  term on the right side of \eqref{eqn:1061}. Putting these together, we have
\begin{equation}\notag
P(A\cap B\cap M^c)\leq \frac{P(M^c)}{\varepsilon P(M)}P(A\cap B\cap M)
\end{equation}
Therefore,
\begin{align}\notag
        P(A\cap B)&=P(A\cap B\cap M^c)+P(A\cap M\cap B)\\\notag
        &\leq \left(\frac{P(M^c)}{\varepsilon P(M)}+1\right)P(A\cap B\cap M)
        \\\notag
        &\overset{\eqref{eqn:1054}}{=}\left(\frac{P(M^c)}{\varepsilon P(M)}+1\right)P(A)P(B)P(M)
        \\\notag
        &=\left(\frac{P(M^c)}{\varepsilon }+P(M)\right)P(A)P(B).
    \end{align}
    Combining the above with \eqref{eqn:1058}, we get
    \begin{equation}\notag
    P(M)\leq\frac{P(A\cap B)}{P(A)P(B)}\leq \left(\frac{P(M^c)}{\varepsilon }+P(M)\right)
    \end{equation}
    Since $P(M)$ approaches 1 as $n\to\infty$, this completes the proof of the proposition.
\end{proof}

\begin{rem}
    We could also construct a counterexample to the statement of Theorem $\ref{thm:MainThm}$ without the i.i.d. assumption where we can show that the walk is not ballistic in direction $e_1$ at all. To do this, take a sequence of values of $\varepsilon_k$ converging to 0, such that $\sum_{k}\varepsilon_k<(1-p)$. Let a site $x$ be red if, for any $k\geq1$, the nearest site $y$ above $x$ such that $U_y<\varepsilon_k$ has in fact $U_y<\varepsilon_k^2$. For red sites, set the probability of stepping in direction $e_1$ to 0, so that it is impossible to go ``through'' a wall. Moreover, the expected size of a new wall is infinite, so that the expected time to get around the wall is infinite, yielding zero speed. However, an environment constructed this way would not have mixing properties as good as the one in Example \ref{exmp:Counterexample}.
\end{rem}

\begin{question}
    In Example \ref{exmp:Counterexample}, is the walk transient to the right? If so, can we improve the example to get a walk that is actually transient to the left?
\end{question}

\section{Concluding Remarks and Further Questions}\label{sec:Qs}

Theorem \ref{thm:MainThm} shows that the asymptotic velocity of the walk exhibits a sort of weak stability under certain perturbations of the distribution on environments. We say ``weak stability'', because rather than showing that the asymptotic velocity is stable, we show that almost-sure bounds on the quenched drift are stable as bounds on the asymptotic velocity when the probability is relaxed slightly away from 1. This is why the continuity result of Corollary \ref{cor:continuity} was only be stated for Dirac masses, where the almost-sure asymptotic velocity and the almost-sure bounds on the quenched drift are one and the same. A desirable extension of Theorem \ref{thm:MainThm} would be to show that the asymptotic velocity itself is stable in general.

\begin{question} Let $\mu_b$ and $\mu_r$ satisfy the assumptions of Theorem \ref{thm:MainThm}. Suppose a RWRE where transition probability vectors at every site are i.i.d. and drawn according to $\mu_b$ has deterministic limiting velocity $v_b$. 
    Let $v_{b,r,p}$ be the asymptotic velocity when sites are drawn according to $\mu_b$ with probability $p$ and $\mu_r$ with probability $1-p$. Can Theorem \ref{thm:MainThm} be strengthened to show that for every $\varepsilon>0$ there exists $p^*\in(0,1)$ such that if $p>p^*$, then $|v_{b,r,p}-v_b|<\varepsilon$?   
\end{question}
This can be re-framed as a question about how much we can weaken the assumption in Corollary \ref{cor:continuity} that $\mu$ is a Dirac mass.

Our proof that the walk can be made to spend an arbitrarily small amount of time at red sites makes no use of the bounds on quenched drift at blue sites. However, \textit{a priori}, the existence of the red sites could change \textit{which} blue sites are visited, and in particular how often blue sites with various drifts are visited. This is why the only bounds on limiting velocity we are currently able to get are bounds that we can read off from the almost-sure bounds on the quenched drift at blue sites. The extension discussed above would require translating our control on the frequency with which the walk visits red sites into a control on the extent to which the visits to red sites affect which blue sites are visited and how often.

Another direction to explore is testing the limits of what perturbations are allowed. We currently insist that $\mu_r$ must have a uniform ellipticity assumption in at least two fixed directions, but impose no conditions beyond these. 
\begin{question}
    Can we relax the uniform ellipticity assumption on $\mu_r$ even further?
\end{question}

It is clear that there must be some constraints on $\mu_r$. If a positive fraction of red sites send the walker in direction $e_i$ with probability 1, and a positive fraction send the walker in direction $-e_i$ with probability 1, then eventually the walk will get stuck between two vertices, almost surely. Moreover, even if constraints are put on $\mu_r$ to prevent a walk from getting stuck permamently---a phenomenon we might refer to as \textit{quenched finite traps}---the asymptotic velocity would still be made 0 by the presence of \textit{annealed finite traps}: finite regions with infinite expected exit time under the annealed measure. For discussions of measures where these appear, see \cite{Tournier2009}, \cite[Theorem 2(i)]{Sabot2013}, \cite[Section 2]{Ramirez&Ribeiro2022}. If $\mu_r$ is a measure of the type described in these sources (or any measure that causes annealed finite traps), then the asymptotic velocity of the walk will become zero as soon as the probability of drawing transition probabilities at a given site from $\mu_r$ becomes positive. Nevertheless, one could explore whether the asymptotic direction---rather than the asymptotic velocity---is stable under such perturbations as these. One could also explore whether imposing conditions that prevent annealed finite traps is enough to cause stability of bounds on the asymptotic velocity.

Finally, we showed in Section \ref{sec:Counterexample} that the i.i.d. assumption---at least on the pre-environment---is an important aspect of Theorem \ref{thm:MainThm}. But can we weaken the i.i.d. assumption at all?

\begin{question}\label{quest:RelaxIID}
    Suppose we relax the i.i.d. assumption so that the pre-environment is i.i.d., but conditioned on the pre-environment, the transition probability vectors are not independent. Does Theorem \ref{thm:MainThm} still hold?
\end{question}

The proof of Claim \ref{claim:909}, which involves the walk ${\bf X}^3$ ``getting around'' an extra red site $y$ in order to recouple with ${\bf X}^2$, does not make any use of independence. The proof of Claim \ref{claim:426}, which compares the probability of a cluster of a given size with the probability of a cluster with one more red site, uses independence in a very fundamental way---this is the part that fails in Example \ref{exmp:Counterexample}---but only independence of the pre-environment. It is therefore reasonable to think Theorem \ref{thm:MainThm} might be true under an i.i.d. assumption on the pre-environment only. However, our proof breaks down under this weaker assumption because it is no longer possible to keep the walks ${\bf X}^2$ and ${\bf X}^3$ coupled whenever they are not at site $y$, and also because a quenched expectation involving the number of times a walk started from a neighbor of $y$ hits $x$ before hitting $y$ is no longer necessarily the same under $\omega_2$ and $\omega_3$ (see Equation \ref{eqn:687}).

\bibliographystyle{plain}
\bibliography{default}

\begin{thebibliography}{10}

\bibitem{Berger&Drewitz&Ramirez2014}
Noam Berger, Alexander Drewitz, and {Alejandro F.} Ram{\'i}rez.
\newblock Effective polynomial ballisticity conditions for random walk in
  random environment.
\newblock {\em Commun. Pure Appl. Math.}, 67(12):1947--1973, December 2014.
\newblock Publisher Copyright: {\textcopyright} 2014 Wiley Periodicals, Inc.
  {\textcopyright} 2014 Wiley Periodicals, Inc.

\bibitem{Bolthausen&Sznitman2002}
Erwin Bolthausen and Alain-Sol Sznitman.
\newblock On the static and dynamic points of view for certain random walks in
  random environment.
\newblock volume~9, pages 345--375. 2002.
\newblock Special issue dedicated to Daniel W. Stroock and Srinivasa S. R.
  Varadhan on the occasion of their 60th birthday.

\bibitem{Bolthausen&Zeitouni2007}
Erwin Bolthausen and Ofer Zeitouni.
\newblock Multiscale analysis of exit distributions for random walks in random
  environments.
\newblock {\em Probab. Theory Related Fields}, 138(3-4):581--645, 2007.

\bibitem{Drewitz&Ramirez2011}
Alexander Drewitz and Alejandro~F. Ram\'irez.
\newblock Ballisticity conditions for random walk in random environment.
\newblock {\em Probab. Theory Relat. Fields}, 150(1):61--75, 06 2011.

\bibitem{Grimmett1999}
G.~Grimmett.
\newblock {\em Percolation}.
\newblock Die Grundlehren der mathematischen Wissenschaften in
  Einzeldarstellungen. Springer, 1999.

\bibitem{Guerra&Ramirez2019}
Enrique Guerra and Alejandro Ram\'irez.
\newblock A proof of sznitman's conjecture about ballistic rwre.
\newblock {\em Communications on Pure and Applied Mathematics}, 11 2019.

\bibitem{Hoggstrom&Jonasson2006}
Olle Häggström and Johan Jonasson.
\newblock Uniqueness and non-uniqueness in percolation theory.
\newblock {\em Probability Surveys}, 3(none), jan 2006.

\bibitem{Kalikow1981}
Steven~A. Kalikow.
\newblock Generalized random walk in a random environment.
\newblock {\em Ann. Probab.}, 9(5):753--768, 10 1981.

\bibitem{Key1984}
Eric~S. Key.
\newblock Recurrence and transience criteria for random walk in a random
  environment.
\newblock {\em Ann. Probab.}, 12(2):529--560, 05 1984.

\bibitem{Komorowski&Krupa2003}
Tomasz Komorowski and Grzegorz Krupa.
\newblock The law of large numbers for ballistic, multi-dimensional random
  walks on random lattices with correlated sites.
\newblock {\em Annales de l'Institut Henri Poincare (B) Probability and
  Statistics}, 39(2):263--285, 2003.

\bibitem{Ramirez&Ribeiro2022}
Alejandro~F. Ram{\'i}rez and Rodrigo Ribeiro.
\newblock {Computable criteria for ballisticity of random walks in elliptic
  random environment}.
\newblock {\em Electronic Journal of Probability}, 27(none):1 -- 38, 2022.

\bibitem{Rassoul-Agha&Seppalainen2009}
Firas Rassoul-Agha and Timo Sepp{\"a}l{\"a}inen.
\newblock {Almost sure functional central limit theorem for ballistic random
  walk in random environment}.
\newblock {\em Annales de l'Institut Henri Poincaré, Probabilités et
  Statistiques}, 45(2):373 -- 420, 2009.

\bibitem{Sabot2004}
Christophe Sabot.
\newblock {Ballistic random walks in random environment at low disorder}.
\newblock {\em The Annals of Probability}, 32(4):2996 -- 3023, 2004.

\bibitem{Sabot2013}
Christophe Sabot.
\newblock Random dirichlet environment viewed from the particle in dimension
  $d\ge3$.
\newblock {\em Ann. Probab.}, 41(2):722--743, 03 2013.

\bibitem{Slonim2021b}
Daniel~J. Slonim.
\newblock Directional transience of random walks in dirichlet environments with
  bounded jumps.
\newblock Preprint, submitted August 2021. https://arxiv.org/abs/2108.11424.

\bibitem{Solomon1975}
Fred Solomon.
\newblock Random walks in a random environment.
\newblock {\em Ann. Probab.}, 3(1):1--31, 02 1975.

\bibitem{Sznitman2001}
Alain-Sol Sznitman.
\newblock On a class of transient random walks in random environment.
\newblock {\em Ann. Probab.}, 29(2):724--765, 04 2001.

\bibitem{Sznitman2002}
Alain-sol Sznitman.
\newblock An effective criterion for ballistic behavior of random walks in
  random environment.
\newblock {\em Probab. Theory Related Fields}, 122, 02 2002.

\bibitem{Sznitman2003}
Alain-Sol Sznitman.
\newblock On new examples of ballistic random walks in random environment.
\newblock {\em The Annals of Probability}, 31(1):285--322, 2003.

\bibitem{Sznitman&Zerner1999}
Alain-Sol Sznitman and Martin Zerner.
\newblock {A Law of Large Numbers for Random Walks in Random Environment}.
\newblock {\em The Annals of Probability}, 27(4):1851 – 1869, 1999.

\bibitem{Tournier2009}
Laurent Tournier.
\newblock Integrability of exit times and ballisticity for random walks in
  dirichlet environment.
\newblock {\em Electron. J. Probab.}, 14:431--451, 2009.

\bibitem{Zeitouni2004}
O.~Zeitouni.
\newblock Random walks in random environment.
\newblock In J.~Picard, editor, {\em Lecture notes in probability theory and
  statistics: École d’été de probabilités de Saint-Flour XXXI-2001},
  volume 1837 of {\em Lect. Notes Math.}, pages 190--313. Springer, 2004.

\bibitem{Zerner2002}
Martin Zerner.
\newblock A non-ballistic law of large numbers for random walks in i.i.d.
  random environment.
\newblock {\em Electron. Commun. Probab.}, 7:191--197, 2002.

\bibitem{Zerner1998}
Martin P.~W. Zerner.
\newblock {Lyapounov exponents and quenched large deviations for
  multidimensional random walk in random environment}.
\newblock {\em The Annals of Probability}, 26(4):1446 -- 1476, 1998.

\bibitem{Zerner&Merkl2001}
Martin P.~W. Zerner and Franz Merkl.
\newblock A zero-one law for planar random walks in random environment.
\newblock {\em Ann. Probab.}, 29(4):1716--1732, 10 2001.

\bibitem{Zerner2000}
Martin~P.W Zerner.
\newblock Velocity and lyapounov exponents of some random walks in random
  environment.
\newblock {\em Annales de l'Institut Henri Poincare (B) Probability and
  Statistics}, 36(6):737--748, 2000.

\end{thebibliography}

\end{document}